\documentclass[10pt,a4paper]{article}
\usepackage{latexsym,fullpage,amsfonts,amssymb,amsmath,amscd,graphics,epic,mathrsfs,theorem}
\usepackage[all]{xy}

\baselineskip 18pt  \theoremstyle{plain}
\newtheorem{lemma}{Lemma}[section]
\newtheorem{proposition}[lemma]{Proposition}
\newtheorem{remark}[lemma]{Remark}
\newtheorem{example}[lemma]{Example}
\newtheorem{theorem}[lemma]{Theorem}
\newtheorem{definition}[lemma]{Definition}
\newtheorem{notation}[lemma]{Notation}
\newtheorem{corollary}[lemma]{Corollary}
\newtheorem{algorithm}[lemma]{Algorithm}

{\theorembodyfont{\rmfamily} 

\makeatletter

\def\Lsection{%
  \@startsection{Lsection}{1}{\z@}{-18\p@ \@plus 2\p@ \@minus 2\p@}%
    {6\p@}{\normalfont\normalsize\scshape}}

\setlength\bibindent{0em}
\renewenvironment{thebibliography}[1]
     {\Lsection*{References}%
      \small\list{\@biblabel{\@arabic\c@enumiv}}%
           {\settowidth\labelwidth{\@biblabel{#1}}%
            \leftmargin\labelwidth
            \setlength\labelsep{4pt}
            \advance\leftmargin\labelsep
            \usecounter{enumiv}%
            \let\p@enumiv\@empty
            \renewcommand\theenumiv{\@arabic\c@enumiv}}%
      \sloppy
      \clubpenalty4000
      \@clubpenalty \clubpenalty
      \widowpenalty4000%
      \sfcode`\.\@m}
     {\def\@noitemerr
       {\@latex@warning{Empty `thebibliography' environment}}%
      \endlist}

\makeatother

\begin{document}
\newcommand{\pperp}{\hbox{$\perp\hskip-6pt\perp$}}
\newcommand{\N}{{\mathbb N}}
\newcommand{\PP}{{\mathbb P}}
\newcommand{\Z}{{\mathbb Z}}
\newcommand{\Q}{{\mathbb Q}}
\newcommand{\R}{{\mathbb R}}
\newcommand{\C}{{\mathbb C}}
\newcommand{\K}{{\mathbb K}}
\newcommand{\F}{{\mathbb F}}
\newcommand{\proofend}{\hfill$\Box$\bigskip}
\newcommand{\eps}{{\varepsilon}}
\newcommand{\ko}{{\mathcal O}}
\newcommand{\wx}{{\widetilde x}}
\newcommand{\wz}{{\widetilde z}}
\newcommand{\wa}{{\widetilde a}}
\newcommand{\wF}{{\widetilde F}}
\newcommand{\wpsi}{{\widetilde \psi}}
\newcommand{\wG}{{\widetilde G}}
\newcommand{\wv}{{\widetilde v}}
\newcommand{\bz}{{\boldsymbol z}}
\newcommand{\bp}{{\boldsymbol p}}
\newcommand{\wy}{{\widetilde y}}
\newcommand{\we}{{\widetilde e}}
\newcommand{\wg}{{\widetilde g}}
\newcommand{\wb}{{\widetilde b}}
\newcommand{\wu}{{\widetilde u}}
\newcommand{\wc}{{\widetilde c}}
\newcommand{\bi}{{\omega}}
\newcommand{\bx}{{\boldsymbol x}}
\newcommand{\Log}{{\operatorname{Log}}}
\newcommand{\pr}{{\operatorname{pr}}}
\newcommand{\Graph}{{\operatorname{Graph}}}
\newcommand{\jet}{{\operatorname{jet}}}
\newcommand{\Tor}{{\operatorname{Tor}}}
\newcommand{\sqh}{{\operatorname{sqh}}}
\newcommand{\const}{{\operatorname{const}}}
\newcommand{\Arc}{{\operatorname{Arc}}}
\newcommand{\Sing}{{\operatorname{Sing}}}
\newcommand{\Span}{{\operatorname{Span}}}
\newcommand{\Aut}{{\operatorname{Aut}}}
\newcommand{\Ker}{{\operatorname{Ker}}}
\newcommand{\Int}{{\operatorname{Int}}}
\newcommand{\Aff}{{\operatorname{Aff}}}
\newcommand{\Area}{{\operatorname{Area}}}
\newcommand{\val}{{\operatorname{Val}}}
\newcommand{\conv}{{\operatorname{conv}}}
\newcommand{\rk}{{\operatorname{rk}}}
\newcommand{\ow}{{\overline w}}
\newcommand{\ov}{{\overline v}}
\newcommand{\ks}{{\cal S}}
\newcommand{\red}{{\operatorname{red}}}
\newcommand{\kc}{{\cal C}}
\newcommand{\ki}{{\cal I}}
\newcommand{\kj}{{\cal J}}
\newcommand{\ke}{{\cal E}}
\newcommand{\kz}{{\cal Z}}
\newcommand{\tet}{{\theta}}
\newcommand{\Del}{{\Delta}}
\newcommand{\bet}{{\beta}}
\newcommand{\mm}{{\mathfrak m}}
\newcommand{\kap}{{\kappa}}
\newcommand{\del}{{\delta}}
\newcommand{\sig}{{\sigma}}
\newcommand{\alp}{{\alpha}}
\newcommand{\Sig}{{\Sigma}}
\newcommand{\Gam}{{\Gamma}}
\newcommand{\gam}{{\gamma}}
\newcommand{\Lam}{{\Lambda}}
\newcommand{\lam}{{\lambda}}
\newcommand{\om}{{\omega}}

\title{Geometry of obstructed equisingular families of algebraic hypersurfaces}
\author{Anna Gourevitch\thanks{ Anna Gourevitch, School of Math. Sciences,
 Tel Aviv University, Ramat Aviv, 69978 Tel Aviv, Israel. Email:
 annabin@post.tau.ac.il.} and Dmitry Gourevitch\thanks {Dmitry Gourevitch, Faculty of Mathematics and Computer
Science, The Weizmann Institute of Science POB 26, Rehovot 76100,
ISRAEL. E-mail: dmitry.gourevitch@weizmann.ac.il. \smallskip
\newline Keywords: T-smoothness, equisingular family, projective
hypersurfaces, quasihomogeneous singularity.
\newline MSC codes: 14J10,14J17,14J70, 32G11.}   }

%
%

%
%
%
%
%
%

\date{}
\maketitle

\begin{abstract}


We study geometric properties of certain obstructed equisingular
families of projective hypersurfaces with quasihomogeneous
singularity with emphasis on smoothness, reducibility, being
reduced, and having expected dimension.

In the case of minimal obstructedness, we give a detailed
description of such families corresponding to quasihomogeneous
singularities.

Next we study the behavior of these properties with respect to
stable equivalence of singularities.

We show that under certain conditions, stabilization of
singularities ensures the existence of a reduced component of
expected dimension. For minimally obstructed families the whole
family becomes irreducible.

As an application we show that if the equisingular family of a
projective hypersurface $H$ has a reduced component of expected
dimension then the deformation of $H$ induced by the equisingular
family $|H|$ is complete with respect to one-parameter
deformations.
\end{abstract}
\setcounter{tocdepth}{1} \tableofcontents

\section{Introduction}
The study of equisingular families of algebraic curves and
hypersurfaces with given invariants and given set of singularities
is an old, but still attractive and widely open problem. Already
at the beginning of the 20th century, the foundation was made in
the works of Pl\"{u}cker, Severi, Segre and Zariski. Later  the
theory of equisingular families has been in focus of the numerous
studies by algebraic geometers and has found important
applications in singularity theory, topology of complex algebraic
curves and surfaces, and in real algebraic geometry.

This paper is devoted to the study of the so-called
\emph{obstructed} families of projective hypersurfaces, of a given
degree, having one isolated singularity of prescribed type.

Let $\Sigma$ be a smooth projective variety over the complex field
$\C$. Let $D$ be an ample divisor on $\Sigma$. Denote by
$V=V_{|D|}(S_{1},\ldots,S_{r})$ the set of hypersurfaces in the
linear system $|D|$ having $r$ singular points of analytic types
$S_{1},\ldots,S_{r}$ (as their only singularities). One knows that
$V_{|D|}(S_{1},\ldots,S_{r})$ can be identified with a (locally
closed) subscheme ("\emph{equisingular stratum"}) in the Hilbert
scheme of hypersurfaces on $\Sig$. The main questions concerning
this space are
\begin{itemize}\item \emph{Existence problem:} Is $V_{|D|}(S_{1},\ldots,S_{r})$ non-empty,
that is, does there exist a hypersurface $F\in|D|$ with the given
collection of singularities?
\item \emph{Smoothness problem:} If $V_{|D|}(S_{1},\ldots,S_{r})$ is non-empty, is it
 smooth?

\item \emph{Dimension problem:} If $V_{|D|}(S_{1},\ldots,S_{r})$ is non-empty, does it
have the "expected" dimension (expressible via local invariants of
the singularities)?

\item \emph{Irreducibility problem:} Is $V_{|D|}(S_{1},\ldots,S_{r})$
irreducible?

\item \emph{Versality problem:} Is the deformation of the
multisingularity of a hypersurface $H \in
V_{|D|}(S_{1},\ldots,S_{r})$ induced by the linear system $|D|$
versal (see Section \ref{SecVerDef})?
\end{itemize}

If $V_{|D|}(S_{1},\ldots,S_{r})$ is non-empty, smooth and has
expected dimension, it is said to be "T-smooth". It is known that
in this case the deformation induced by $|D|$ is versal.

The case of plane nodal curves has been settled completely. In
1920, Severi gave answers to the first three questions: there
exists an irreducible plane curve of degree $d$ having $n$ nodes
as their only singularities if and only if
$$ 0 \leq n \leq \frac{(d-1)(d-2)}{2}.$$ Furthermore, if
$V^{irr}_d(n \cdot A_1)$ is non-empty, then it is smooth of the
expected dimension $\frac{d(d+3)}{2}-n$. In 1985, Harris proved
that $V^{irr}_d(n \cdot A_1)$ is irreducible.

Already in the case of curves with nodes and cusps there is no
complete answer. Segre \cite{Seg,Tan} gave an example of an
equisingular stratum of such curves (having $6m^2$ cusps as their
only singularities) which has a component of non-expected
dimension. Zariski \cite{Zar} gave the first example of a
reducible equisingular stratum (of sextic curves with 6 cusps).
Finally, Wahl \cite{Wah} gave the first example of a non-smooth
equisingular stratum (of curves of degree 104 having 3636 nodes
and 900 cusps). Some other examples can be found in \cite{Lue},
\cite{GLS2}, \cite{Kei}, \cite{Mar}.

So far, the main effort in the study of equisingular families has
been concentrated on obtaining criteria for the T-smoothness. In
turn, the obstructed (i.e., non-T-smooth) equisingular families
and non-versal deformations have not been studied systematically.

The deformation theory leads to the following result: Suppose that
$h^1({\mathcal O}_\Sigma(D))=0$. Then the variety
$V_{|D|}(S_{1},\ldots,S_{r})$ is T-smooth at $H\in
V_{|D|}(S_{1},\ldots,S_{r})$ if and only if
$$h^{1}(\kj_{Z^{ea}(H)/\Sigma}(H))=0,$$
where $Z^{ea}(H)$ is a certain zero-dimensional scheme and
$\kj_{Z^{ea}(H)/\Sigma}(H)$ is its defining ideal (see precise
definition in Section \ref{Zero-dim}). In our research we focus on
the minimal obstructedness case, i.e.
$h^{1}(\kj_{Z^{ea}(H)/\Sigma}(H))=1.$

One of the interesting recent examples is due to du Plessis and
Wall \cite{DPW}:
\begin{example}\label{e2}(a) For any $d\geq5$ the curve
$C\subset\PP^2$ given by the equation
$(x_1^d+x_2^5x_0^{d-5}+x_2^d=0)$ has a unique singular point
$z=(0:0:1)$ with Tjurina number $\tau(C,z)=4d-4$, and satisfies
$$h^{1}(\kj_{Z^{ea}(C)/\PP^2}(d))>0$$
(b) Denote by $S$ the analytic type of the plane curve singularity
$(C,z)$ in (a). If $d\geq10$ then the family  $V_{d}(S)$ is
singular at $C$.
\end{example}


In \cite{Gou1} this example has been generalized and studied in
detail. The following result has been obtained:

\begin{example}\label{Gur1} Let $C$ be a projective plane curve, given
in local coordinates $x=\frac{x_1}{x_0}$, $y=\frac{x_2}{x_0}$ by
the  equation $x^{k}+y^{l}=0$. Suppose for convenience $k\geq l$.
Let $d=k+l-5$ and let $V_{d,C}(S)$ be the germ  at $C$ of the
equianalytic family of plane curves $V_{d}(S)$, where $S$ is the
analytic type of the plane curve singularity $(C,z)$. For any $k,
l\geq 5$ such that $d>5$, $V_{d,C}(S)$ is non-T-smooth and
$h^1(\kj_{Z^{ea}(C)/\mathbb{P}^{2}}(d))=1$.\\ Furthermore

(i) If $d=6$ (i.e. $l=5$, $k=6$), the germ $V_{6,C}(S)$  is
non-reduced. It is a double $[V_{6,C}(S)]_{red}$, and
$[V_{6,C}(S)]_{red}$ is smooth of expected codimension.

 (ii) If $d=7$ (i.e. $l=5$, $k=7$ or $k=l=6$), the germ $V_{7,C}(S)$ is
reducible and  decomposes into two smooth components of expected
codimension that intersect non-transversally with multiplicity
one. The intersection locus is smooth. Moreover, the sectional
singularity is of type $A_{1}$.

 (iii) If $d\geq 8$, the germ $V_{d,C}(S)$ is
a reduced irreducible non-smooth variety of expected codimension
which has a smooth singular locus  with sectional singularity of
type $A_{1}$.

\end{example}

\subsection{Main results and methods}

The first result of this paper is a generalization of Example
\ref{Gur1} to quasihomogeneous hypersurface singularities in
$\PP^n$. In particular we have obtained a new example of a smooth
equisingular family of non-expected dimension: $V_{3,H}(S)$ where
$H$ is given by the local equation $x^3+y^3+z^3+w^3=0$. For
precise formulation see Theorem \ref{QHomN}.

The next question that naturally arose was the behavior of these
geometric properties of equianalytic families, with respect to the
stabilization of the singularities. We found out that though
stabilization preserves both the Tjurina algebra and
$h^{1}(\kj_{Z^{ea}(H)/\PP^n}(H))$, it can change the geometry of
the equianalytic family radically. We have shown that for any
hypersurface $H$ of degree $d$ satisfying
$h^{1}(\kj_{Z^{ea}(H)/\PP^n}(2d-2))=0$, after adding enough
squares the obtained family has a reduced component of expected
dimension. We have also shown that the condition
$h^{1}(\kj_{Z^{ea}(H)/\PP^n}(2d-2))=0$ always holds for plane
curves. For precise formulation see Theorem \ref{Squares}.

The next result is concerned with deformation theory. Suppose that
the germ of the equianalytic family has a reduced irreducible
component of expected dimension. In this case the family is
T-smooth at every of its regular points which lies in that
component. It means that our singular hypersurface $H$ has a
deformation $\mathfrak{X}\to T$ such that for $t\neq t_0$ the
deformation of $\mathfrak{X}_t$ induced by the linear system $|H|$
is versal. We show that this implies that the deformation of $H$
induced by the linear system $|H|$ is 1-complete (see Section
\ref{DefTheory} for precise definition). For precise formulation
see Theorem \ref{DefMeaning}.

The methods that we use are the technique of cohomologies of ideal
sheaves of zero-dimensional schemes associated with analytic types
of singularities, methods for their calculation, and
$H^{1}$-vanishing theorems. We also use the algorithms of computer
algebra (see \cite{GP}) as a technical tool in the proof of the
theorems.

\subsection*{The structure of the paper}

This paper is organized in the following way. Section 2 is
dedicated to the formulation of the necessary notions and
background.

In Subsection \ref{SecSingNot} we introduce the notions of
singularity theory such as \textbf{analytic singularity types},
\textbf{quasihomogeneous} singularities, \textbf{zero-dimensional
schemes associated with singularities} and the \textbf{Castelnuovo
function}. The main theorems of this section are the Mather-Yau
theorem (Theorem \ref{M-Y}), finite determinacy theorem (Theorem
\ref{F-D}), Theorem \ref{Saito} on quasihomogeneous singularities,
Theorem \ref{Arn} on semiquasihomogeneous polynomials, Theorem
\ref{p3} which gives cohomological criteria of T-smoothness, and
Lemmas \ref{Cast} and \ref{Davis} on the Castelnuovo function.

In Subsection \ref{DefTheory} we introduce the notions of
deformation theory such as \textbf{complete} and \textbf{versal}
deformations. The most important statement for us in this section
is Corollary \ref{VerDef} on versality of the deformation induced
by a complete linear system.

In Subsection \ref{SecComp} we introduce notions and algorithms of
computer algebra. The most important notions for us are the notion
of \textbf{normal form} and the \textbf{\textsc{RedNFBuchberger}}
algorithm for its computation.

In Section \ref{Main} we formulate our \textbf{main results}. The
first result deals with equianalytic families of hypersurfaces
with quasihomogeneous singularities of minimal obstructedness. The
second result is about stable properties of obstructed
equianalytic families. In Subsection \ref{SecOurDef} we give an
application of the obtained results to the deformation theory.

In Section \ref{ProofQHomN} we prove the theorem on families of
hypersurfaces with quasihomogeneous singularities of minimal
obstructedness.

In Section \ref{ProofSquares} we prove the theorem on stable
properties of obstructed equianalytic families.

\subsection*{Acknowledgements}

First of all we would like to thank  \textbf{Prof. Eugenii
Shustin} for  supervising this work. 

We would also like to thank \textbf{Prof. Joseph Bernstein} for
teaching us basics of Algebraic Geometry, \textbf{Prof.
Gert-Martin Greuel} for fruitful discussions and \textbf{Dr.
Dmitry Kerner} for useful remarks.

Finally, we wish to thank the referee for very careful reading and
useful remarks.

The first-named author was supported by Hermann-Minkowski-Minerva
Center for Geometry at the Tel-Aviv University, by the grant no.
465/04 from the Israel Science Foundation. The second-named author
was supported by GIF grant no. 861/05 and ISF grant no. 1438/06.
%
%
\section{Preliminaries and notations}
\subsection{Notions of singularity theory} \label{SecSingNot}

In this section we describe  the types of isolated hypersurface
singularities considered throughout the paper.
\subsubsection{Analytic types of hypersurface singularities}

\begin{definition}
{\rm  Let $\Sigma$ be an n-dimensional smooth projective variety.
Two germs $(F,z)\subset (\Sigma,z)$ and $(G,w)\subset(\Sigma,w)$
of isolated hypersurface singularities are said to be
\emph{analytically equivalent} if there exists a local analytic
isomorphism $(\Sigma,z)\rightarrow(\Sigma,w)$ mapping $(F,z)$ to
$(G,w)$. The corresponding equivalence classes are called
\emph{analytic types}.}
\end{definition}

\begin{notation}
{\rm We denote by $\mathbb{C}\{x_1,\ldots,x_n\}$, or in short
$\mathbb{C}\{x\}$, the algebra of convergent power series in n
variables.}
\end{notation}
\begin{definition}
{\rm  Series $f,g\in \mathbb{C}\{x\}$ are said to be contact
equivalent if there exists an automorphism $\phi$ of
$\mathbb{C}\{x\}$ and a unit $u\in\mathbb{C}\{x\}^*$ such that
$f=u\cdot \phi(g)$. We denote $f\overset{c}{\sim} g$.}
\end{definition}
Note that polynomials $f$ and $g$ are contact equivalent if and
only if the corresponding germs $(f^{-1}(0),0)$ and
$(g^{-1}(0),0)$ are analytically equivalent.
\begin{definition}\label{r1}{\rm  Let $S$ be an analytic type of reduced hypersurface
 singularities represented by $(H,z)\subset (\Sigma,z)$ and $f\in
\mathbb{C}\{x_1,\ldots,x_n\}$ be a local equation for $(H,z)$.
Define the jacobian of $f$ by $j(f)=\langle \frac{\partial
f}{\partial x_{1}},\ldots,\frac{\partial f}{\partial
x_{n}}\rangle.$ The analytic algebras $$
M_{f}:={\mathbb{C}\{x\}}/j(f),\quad
T_{f}:={\mathbb{C}\{x\}}/\langle f,j(f)\rangle$$ are called the
\emph{Milnor } and \emph{Tjurina} algebra of $f$, respectively,
and the numbers $$\mu(S)=\mu(H,z):=\dim_{\mathbb{C}}M_{f},\quad
\tau(S)=\tau(H,z):=\dim_{\mathbb{C}}T_{f}$$ are called the
\emph{Milnor } and \emph{Tjurina} numbers of $S$, respectively.}

\end{definition}

The following theorem shows that the Tjurina algebra is a complete
invariant of an analytic singularity type.

\begin{theorem}[Mather-Yau] \label{M-Y} Let
$f,g\in m \subset \mathbb{C}\{x\}$. The following are
equivalent:\\(a) $f\overset{c}{\sim} g$;\\(b) for all $b\geq 0$,
$\mathbb{C}\{x\}/\langle f,m^bj(f)\rangle \cong
\mathbb{C}\{x\}/\langle g,m^bj(g)\rangle$ as
$\mathbb{C}$-algebras;\\
(c) there is some $b\geq 0$ such that $\mathbb{C}\{x\}/\langle
f,m^bj(f)\rangle \cong \mathbb{C}\{x\}/\langle g,m^bj(g)\rangle$
as $\mathbb{C}$-algebras. \\In particular, $f\overset{c}{\sim} g$
iff $T_f\cong T_g$.
\end{theorem}
\emph{Proof.} See \cite{MaY} for the case of an isolated
singularity and $b=0,1$ or \cite{GLS}, Theorem 2.26 for the
general case.

%
%
\subsubsection{Finite determinacy} The aim of this section is to
show that an isolated hypersurface singularity is already
determined by its Taylor series expansion up to a sufficiently
high order.
\begin{definition}\label{jet}{\rm
For $f\in \mathbb{C}\{x\}$ we define the \emph{$k-$jet} of $f$ by
$$jet(f,k):=f^{(k)}:=\text{image of $f$ in }\mathbb{C}\{x\}/m^{k+1}.$$
We identify $f^{(k)}$ with the power series expansion of $f$ up to
(and including) order $k$.}

\end{definition}
\begin{definition}\label{DFD}{\rm
$f\in \mathbb{C}\{x\}$ is called \emph{contact $k$-determined} if
for each $g\in \mathbb{C}\{x\}$  with $f^{(k)}=g^{(k)}$ we have
$f\overset{c}{\sim} g$. The minimal such $k$ is called
\emph{contact determinacy of $f$}.}

\end{definition}
\begin{theorem}\label{F-D}(Finite determinacy theorem) Let
$f\in m \subset \mathbb{C}\{x\}$.\\ $f$ is contact $k-$determined
if $m^{k+1}\subset m^2\cdot j(f)+m \langle f\rangle.$
\end{theorem}
\emph{Proof.} This theorem is well known. See, for instance,
\cite{GLS}, Theorem 2.23.
\begin{corollary}\label{CFD} If
$f\in m \subset \mathbb{C}\{x\}$, has an isolated singularity with
Tjurina number $\tau$, then $f$ is contact $\tau +1-$determined .
\end{corollary}
\subsubsection{Stable equivalence}\label{SecStEqui}
\begin{definition}{\rm Let $f\in \mathbb{C}\{x_{1},\ldots,x_{l}\}$
and $g\in \mathbb{C}\{x_{1},\ldots,x_{k}\}$. We say that $f$ is
\emph{stably contact equivalent} to $g$ if they become contact
equivalent after addition with non-degenerate quadratic forms of
additional variables. In other words,
$$f(x_{1},\ldots,x_{l})+x_{l+1}^2+\ldots+x_{n}^2 \overset{c}{\sim} g(x_{1},\ldots,x_{k})+x_{k+1}^2+\ldots+x_{n}^2$$  }
\end{definition}

\begin{theorem}
Polynomials of the same number of variables are stably contact
equivalent if and only if they are contact equivalent.
\end{theorem}
\emph{Proof.} See, for instance, \cite{AGV}, chapter
II section 11. \\
We will use a more general lemma:
\begin{lemma} \label{MathYauCor}
(1) Let $f \in \C\{x_1,...,x_n\}$. Let $g=f+x_{n+1}^2(1+h)$ where
$h \in m \subset \C\{x_1,...,x_{n+1}\}$. Then the Tjurina algebra
 $T_f$ of $f$ is isomorphic to the Tjurina algebra $T_g$ of $g$.\\
(2) Let $f_1,f_2 \in \C\{x_1,...,x_n\}$. Let
$g_i=f_i+x_{n+1}^2(1+h_i) \, , i=1,2$ where $h_i \in m \subset
\C\{x_1,...,x_{n+1}\}$. Then $f_1 \overset{c}{\sim} f_2$ if and
only if $g_1 \overset{c}{\sim} g_2$
\end{lemma}
\emph{Proof.} (1) In $T_g$, $$0=\frac{\partial g}{\partial
x_{n+1}}=x_{n+1}(2+2h+x_{n+1}\frac{\partial h}{\partial
x_{n+1}}).$$ Since $2+2h+x_{n+1}\frac{\partial h}{\partial
x_{n+1}}$ is invertible in $T_g$, the latter implies $x_{n+1}=0$.
Hence
$T_f$ is isomorphic to $T_g$.\\
(2) Follows from (1) and from the Mather-Yau theorem (\ref{M-Y}).
\proofend

\subsubsection{Quasihomogeneous singularities}
\begin{definition}\label{10} {\rm Let
$f=\sum_{\alpha\in\mathbb{Z}^{n}_{\geq 0}}a_{\alpha}x^{\alpha}\in\mathbb{C}[x_{1},\ldots,x_{n}].$

(i) The polynomial $f$ is called \emph{quasihomogeneous of type}
$$(w;d)=(w_{1},\ldots,w_{n};d)$$ if $w_{i}$, $d$ are positive
integers satisfying $$\langle
w,\alpha\rangle=w_{1}\alpha_{1}+\ldots +w_{n}\alpha_{n}=d$$ for
each $\alpha\in\mathbb{Z}^{n}_{\geq 0}$ with $a_{\alpha}\neq 0$.

(ii) An isolated hypersurface singularity $(H,x)\subset(
\mathbb{C}^{n},x)$ is called \emph{quasihomogeneous} if there
exists a quasihomogeneous polynomial $f\in
\mathbb{C}[x_1,\ldots,x_n]$ such that ${\cal
O}_{H,x}\cong\mathbb{C}\{x\}/\langle f\rangle$.}
\end{definition}
\begin{example}\label{CanQHomDef} {\rm
Let  $f=\sum \limits _{i=1}^{n}x_i^{\alpha_i}$. Then $f$ is a
quasihomogeneous polynomial of type $(1/\alp_1,\dots,1/\alp_n;1)$.
Such polynomials are called \emph{canonical quasihomogeneous}.}
\end{example}

\begin{lemma}
Let $f\in \C[x_1,...,x_n]$ be quasihomogeneous and $g\in
\C\{x_1,...,x_n\}$ be arbitrary. Then $f$ and $g$ are contact
equivalent if and only if there exists an analytic local
diffeomorphism $\phi$ that maps $g$ to $f$.
\end{lemma}
\emph{Proof}. Let $f$ be quasihomogeneous of type
$((w_{1},\ldots,w_{n});d)$. If $\overset{c}{\sim} g$ then there
exists a unit $u \in \C\{x\}^*$ and an automorphism $\psi \in
Aut\C\{x\}$ such that $u \cdot f = \psi(g)$. Choose a $d$-th root
$u^{1/d} \in \C\{x\}$. Now we take
$$ \phi:\C\{x\} \to \C\{x\}, \quad x_i \mapsto u^{w_i/d} \cdot
x_i. $$ \proofend

A quasihomogeneous polynomial $f$ of type $(w;d)$ obviously
satisfies the relation $$ d\cdot
f=\sum_{i=1}^{n}w_{i}x_{i}\frac{\partial f}{\partial x_{i}},$$
that implies that $f$ is contained in $j(f)$, hence, for
quasihomogeneous isolated hypersurface singularities $\mu=\tau$.
For an isolated singularity the converse  also holds. More
precisely, K.Saito proved the following theorem.

\begin{theorem}[\cite{Sai}]\label{Saito}
Let $f\in\mathbb{C}\{x_{1},\ldots,x_{n}\}$ and suppose that $f \in
j(f)$. Then there exists a quasihomogeneous polynomial $g$ 
and an analytic local diffeomorphism $\phi$ that maps $g$ to $f$.
Moreover, the normalized quasihomogeneity type $\frac{w}{d}$ of
$g$ is defined uniquely up to permutation.
\end{theorem}
We will use the following important theorem:
\begin{theorem} \label{Arn}
Let $f$ be a quasihomogeneous polynomial of weighted degree $d$
and $g$ be a polynomial such that the weighted degrees of  all its
terms are greater than $d$. Let $e_1,e_2,...,e_{\mu}$ be a system
of monomials that forms a basis of the Milnor algebra $M_f$ of
$f$. Then $e_1,e_2,...,e_{\mu}$ form a basis of $M_{f+g}$ as well.
In particular, $\mu(f+g)=\mu(f)$.
\end{theorem}
For a proof see \cite{AGV}, chapter II Section 12.2.

\begin{corollary} \label{Var}
Let $H$ be a projective hypersurface defined by a quasihomogeneous
polynomial $f$ of weighted degree $d$. Let $V_H$ denote the germ
at $H$ of the equisingular family of $H$. Let $H_F \in V_H$ be a
hypersurface defined by a polynomial $F$. Suppose that all the
terms of $F - f$ of weighted degree less than or equal to $d$ are
elements of a monomial basis of the Milnor algebra of $f$.

Then $F - f$ has no terms of weighted degree less than or equal to
$d$.
\end{corollary}
\emph{Proof}. Decompose $F = f + g + \sum_{i=1}^{r} \lambda_i
e_i$, where $e_i$ are all the elements of a monomial basis of the
Milnor algebra $M_f$ of weighted degree less than or equal to $d$,
and $g$ has no terms which have weighted degree less than or equal
to $d$. By the theorem, $M_{f+g}$ and $M_f$ have the same basis,
and $f$ belongs to the $\mu$-constant stratum of $f+g$. Hence
$F=(f+g)+ \sum_{i=1}^{r} \lambda_i e_i$ belongs to the
$\mu$-constant stratum of $f+g$. On the other hand, by \cite{Var},
the affine space $(f+g)+ \Span\{e_i\}_{i=1}^r$ is transversal to
the $\mu$-constant stratum of $f+g$ and hence $\sum_{i=1}^{r}
\lambda_i e_i =0$. \proofend

\subsubsection{Newton polytope}
\begin{definition}{\rm
Let $f=\sum_{\alp\in \mathbb{Z}_{\geq 0}^n}a_{\alp}x^{\alp}\in
\mathbb{C}\{x\}$, $a_0=0$. The convex hull in $\mathbb{R}^n$ of
the support of $f$,
$$ \Delta(f):=\conv \{\alp \in \mathbb{Z}_{\geq 0}^n \mid a_{\alp}\neq 0\},$$
 is called \emph{the Newton polytope of $f$}.}
\end{definition}
If $f$ is quasihomogeneous with weight $w$ of degree $d$ then the
Newton polytope $\Delta(f)$ will lie in the hyperplane $\{\alp \in
\mathbb{Z}_{\geq 0}^n \mid \sum\alp_iw_i=d\}$. In particular, for
a quasihomogeneous polynomial of two variables its Newton polygon
is a line segment.

\subsubsection{Zero-dimensional schemes associated with singularities}\label{Zero-dim}

\begin{definition}{\rm Let $\Sig$ be a smooth projective variety and $H\subset\Sig$
a reduced hypersurface with singular locus
$Sing(H)=\{z_{1},...,z_{r}\}$. We define $Z^{ea}(H)\subset\Sigma$
 to be the zero-dimensional scheme, concentrated at $Sing(H)$, given by
the Tjurina ideals
$$\mathcal{J}_{Z^{ea}(H)/\Sigma,z_{i}}=\langle f_i,\frac{\partial
f_i}{\partial x_1},\dots , \frac{\partial f_i}{\partial
x_n}\rangle \subset\mathcal{O}_{\Sigma,z_{i}},$$ where
$f_i(x_1,\dots ,x_n)$ is a local equation of $H$ in a neighborhood
of $z_i$. We denote by $\mathcal{J}_{Z^{ea}(H)/
\Sigma}\subset\mathcal{O}_{\Sigma}$ the corresponding ideal
sheaf.}
\end{definition}
The degree of $Z^{ea}(H)$ is
$$\textrm{deg}\,Z^{ea}(H)=\sum_{z_i\in
Sing(H)}dim_{\mathbb{C}}\mathcal{O}_{\Sigma,z_{i}}/\mathcal{J}_{Z^{ea}(H)/\Sigma,z_{i}}.$$
\begin{theorem}\label{p3} Let $H\subset \mathbb{P}^{n}$ be a
reduced hypersurface of degree $d$ with precisely $r$
singularities $z_{1},...,z_{r}$ of analytic 
types $S_{1},...,S_{r}$. \\(a)
$H^{0}(\mathcal{J}_{Z^{ea}(H)/\mathbb{P}^{n}}(d))/H^{0}(\mathcal{O}_{\mathbb{P}^{n}})$
is isomorphic to the Zariski tangent space to
$V_{|H|}(S_{1},...,S_{r})$ at $H$. Here,
$H^{0}(\mathcal{O}_{\mathbb{P}^{n}})$ is embedded into
$H^{0}(\mathcal{J}_{Z^{ea}(H)/\mathbb{P}^{n}}(d))$ via
multiplication by the equation of $H$.\\
(b)
$h^{0}(\mathcal{J}_{Z^{ea}(H)/\mathbb{P}^{n}}(d))-h^{1}(\mathcal{J}_{Z^{ea}(H)/\mathbb{P}^{n}}(d))-1\leq
dim(V_{|H|}(S_{1},...,S_{r}),H)\leq
h^{0}(\mathcal{J}_{Z^{ea}(H)/\mathbb{P}^{n}}(d))-1.$\\
(c) $H^{1}(\mathcal{J}_{Z^{ea}(H)/\mathbb{P}^{n}}(d))=0$ if and
only if $V_{|H|}(S_{1},...,S_{r})$ is T-smooth at $H$, that is,
smooth of the expected dimension $\binom{d+n}{n}-1-deg\,
Z^{ea}(H).$

\end{theorem}
\emph{Proof.} See \cite{GrK}.
\begin{lemma}\label{prelim} Let $H\subset \PP^n$ be a hypersurface
and $Z\subset H$ be a zero-dimensional subscheme. Then
$H^{1}(\mathcal{J}_{Z/\mathbb{P}^{n}}(H))\cong
H^{1}(\mathcal{J}_{Z/H}(H))$.
\end{lemma}
\emph{Proof.} The lemma follows from the exact sequence of sheaves
$$0\to \mathcal{J}_{H/\mathbb{P}^{n}}(H)\to
\mathcal{J}_{Z/\mathbb{P}^{n}}(H)\to
\mathcal{J}_{Z/\mathbb{P}^{n}}(H)\otimes \mathcal{O}_H\to 0,$$
from the fact that $\mathcal{J}_{H/\mathbb{P}^{n}}(H)\cong
\mathcal{O}_{\PP^n}$, which implies
$H^1(\mathcal{J}_{H/\mathbb{P}^{n}}(H))=H^2(\mathcal{J}_{H/\mathbb{P}^{n}}(H))=0$,
and from $\mathcal{J}_{Z/\mathbb{P}^{n}}(H)\otimes
\mathcal{O}_H\cong \mathcal{J}_{Z/H}(H)$.

\subsubsection{The Castelnuovo function of a zero-dimensional scheme in $\PP^n$}\label{CastFun}
Let $X\subset \PP^n$ be a zero-dimensional scheme and
$\mathcal{J}_{X/\mathbb{P}^{n}}\subset \mathcal{O}_{\PP^n}$ the
corresponding ideal sheaf.
\begin{definition} {\rm The \emph{Castelnuovo function} of $X$ is
defined as
\begin{align}\mathcal{C}_X:\mathbb{Z}_{\geq 0}&\longrightarrow \mathbb{Z}_{\geq 0} \nonumber \\
d&\longmapsto
h^1(\mathcal{J}_{X/\mathbb{P}^{n}}(d-1))-h^1(\mathcal{J}_{X/\mathbb{P}^{n}}(d)).\nonumber
\end{align}}
\end{definition}
\begin{remark}{\rm Let $X\subset \PP^n$ be a zero-dimensional scheme
and $H\subset \PP^n$ be a generic hyperplane not passing through
the support of $X$. Then we have an exact reduction sequence
$$0\longrightarrow
\mathcal{J}_{X/\mathbb{P}^{n}}(d-1)\longrightarrow
\mathcal{J}_{X/\mathbb{P}^{n}}(d)\longrightarrow
\mathcal{O}_{H}(d)\longrightarrow 0,$$ respectively the
corresponding exact cohomology sequence
$$H^0(\mathcal{J}_{X/\mathbb{P}^{n}}(d-1))\overset{\pi_H}\longrightarrow
 H^0(\mathcal{O}_{H}(d))\longrightarrow
H^1(\mathcal{J}_{X/\mathbb{P}^{n}}(d-1))\longrightarrow
H^1(\mathcal{J}_{X/\mathbb{P}^{n}}(d))\longrightarrow 0.$$ In
particular,
$$\mathcal{C}_X(d)=h^0(\mathcal{O}_{H}(d))-\dim_{\mathbb{C}}\pi_H(H^0(\mathcal{J}_{X/\mathbb{P}^{n}}(d))).$$}
\end{remark}
We associate to $X$ the numbers
\begin{align}
& a(X)=\min\{d\in\mathbb{Z}|\,h^0(\mathcal{J}_{X/\mathbb{P}^{n}}(d)>0)\} \nonumber\\
& b(X)=\min\{d\in
\mathbb{Z}|\,\PP(H^0(\mathcal{J}_{X/\mathbb{P}^{n}}(d))) \text{
has no fixed
component}\} \nonumber\\
& t(X)=\min\{d\in
\mathbb{Z}|\,H^1(\mathcal{J}_{X/\mathbb{P}^{n}}(d))=0)\}.
\nonumber
\end{align}
Here, a fixed component is a divisor $D$ such that every element
of the linear system $|H^0(\mathcal{J}_{X/\mathbb{P}^{n}}(d))|$
contains $D$ as a component. We call the maximal divisor
satisfying this property \emph{the fixed component of}
$|H^0(\mathcal{J}_{X/\mathbb{P}^{n}}(d))|$. The following lemma
contains some basic properties of the Castelnuovo function.
\begin{lemma}\label{Cast} Let $X\subset \PP^n$ be a zero-dimensional scheme
and $H\subset \PP^n$ be a generic hyperplane not passing through
the support of $X$. Then

(a) $\mathcal{C}_X(d)\geq 0$ for all $d$, and $\mathcal{C}_X(d)=
0$ for $d\gg 0$.

(b)  $\mathcal{C}_X(d)\leq h^0(\mathcal{O}_{H}(d))$, with equality
if and only if $d<a(X)$.

(c) $a(X)\leq b(X)\leq t(X)+1$.

(d) $\mathcal{C}_X(d)= 0$ if and only if $d\geq t(X)+1$.

(e) If $Y\subseteq X$ then $\mathcal{C}_Y(d)\leq
\mathcal{C}_X(d)$.

(f)$\mathcal{C}_X(0)+\mathcal{C}_X(1)+\dots+\mathcal{C}_X(d)=\textmd{deg}X-h^1(\mathcal{J}_{X/\mathbb{P}^{n}}(d))$.
\end{lemma}
\emph{Proof.} See \cite{Dav} for curves or \cite{West} for
hypersurfaces.

\begin{lemma}\label{Davis}
Let $Z=C_d \cap C_k$ be the intersection of two plane curves
$C_d,C_k$ of degrees $d$ and $k$ without common components.
Suppose $k\leq d$. Then $$\mathcal{C}_Z(i) \leq k \text{ for
}i\geq 0, \text{ and }\mathcal{C}_Z(d+k-j)=j-1 \text{ for
}j=1,\dots,k+1.$$
\end{lemma}
\emph{Proof.} See e.g. \cite{GLS1}, Lemma 5.4.
%

\subsection{Some notions of the local deformation theory}\label{DefTheory}

\begin{definition}
{\rm Let $(X,x)$ and $(S,s)$ be complex space germs. \emph{A
deformation of $(X,x)$ over $(S,s)$} consists of a flat morphism
$\phi:(\mathscr{X},x) \to (S,s)$ of complex space germs together
with an isomorphism from $(X,x)$ to the fiber of $\phi$, $(X,x)\to
(\mathscr{X}_s,x):=(\phi^{-1}(s),x)$.

$(\mathscr{X},x)$ is called the \emph{total space}, $(S,s)$ the
\emph{base space}, and $(\mathscr{X}_s,x) \cong (X,x)$ the
\emph{special fiber} of the deformation. We denote a deformation
by $$(i,\phi):(X,x)\overset{i}{\hookrightarrow}
(\mathscr{X},x)\overset{\phi}{\rightarrow}(S,s),$$ or simply by
$(\mathscr{X},x)\overset{\phi}{\rightarrow}(S,s)$.
 }
\end{definition}
\begin{definition}
{\rm Let $(i,\phi):(X,x)\overset{i}{\hookrightarrow}
(\mathscr{X},x)\overset{\phi}{\rightarrow}(S,s)$ and
$(i',\phi'):(X,x)\overset{i}{\hookrightarrow}
(\mathscr{X}',x')\overset{\phi'}{\rightarrow}(S',s')$ be two
deformations of $(X,x)$. A morphism of deformations from
$(i,\phi)$ to $(i',\phi')$ consists of two morphisms
$(\psi,\varphi)$ such that the following diagram is commutative.

$$\quad (X,x)\quad$$
 $$i'\swarrow \quad \quad \searrow i $$
$$(\mathscr{X}',x')\quad \stackrel {\psi}{\longrightarrow}\quad (\mathscr{X},x) $$
 $$\phi'\downarrow \quad \quad \quad \quad \quad \quad \downarrow \phi $$
$$(S',s')\quad \stackrel {\varphi}{\longrightarrow}\quad (S,s) $$
Two deformations over the same base are isomorphic if there exists
a morphism $(\psi,\varphi)$ with $\psi$ an isomorphism and
$\varphi$ the identity map. }
\end{definition}

\begin{definition}
{\rm Let $(i,\phi):(X,x)\overset{i}{\hookrightarrow}
(\mathscr{X},x)\overset{\phi}{\rightarrow}(S,s)$  be a deformation
of $(X,x)$ and $\varphi:(T,t)\to (S,s)$ be a morphism of germs.
Denote by $\varphi^*(\mathscr{X},x)$ the fiber product
$(\mathscr{X},x)\times_{(S,s)}(T,t)$. We call
$$\varphi^*(i,\phi):=(\varphi^*i,\varphi^*\phi):(X,x)\overset{\varphi^*i}{\hookrightarrow}
\varphi^*(\mathscr{X},x)\overset{\varphi^*\phi}{\rightarrow}(T,t)$$
the \emph{deformation induced from $(i,\phi)$ by $\varphi$}, or
just \emph{pull-back}; $\varphi$ is called the \emph{base change
map}.
 }
\end{definition}
\begin{proposition}
Let $(X,0)\subset(\mathbb{C}^n,0)$ be a closed subgerm. Then any
deformation $(i,\phi):(X,x)\overset{i}{\hookrightarrow}
(\mathscr{X},x)\overset{\phi}{\rightarrow}(S,s)$ can be embedded.
In other words, there exists a Cartesian diagram
$$(X,0)\quad \stackrel {i}{\hookrightarrow}\quad (\mathscr{X},x) $$
 $$\downarrow \quad \quad \quad \quad \quad \quad \downarrow J $$
$$\quad \quad\quad(\mathbb{C}^n,0)\quad \stackrel {j}{\hookrightarrow}\quad
(\mathbb{C}^n,0)\times(S,s) $$
$$\downarrow \quad \quad \quad \quad \quad \quad \downarrow p $$
$$\{s\}\quad \quad \stackrel {\hookrightarrow}\quad\quad (S,s) $$
where $J$ is a closed embedding, $p$ is the second projection, $j$
is the first inclusion and $\phi = p \circ J$.
%
In particular, the embedding dimension is semicontinuous under
deformations, that is,
\emph{edim}$\phi^{-1}(\phi(y)),y)\leq$\emph{edim}$(X,0)$, for all
$y$ in $\mathscr{X}$ sufficiently close to $x$.
\end{proposition}
\emph{Proof.} See \cite{GLS}, Corollary II.1.6.

\subsubsection{Versal and complete deformations}\label{SecVerDef}  A versal deformation of a
complex space germ is a deformation which contains basically all
information about any possible deformation of this germ. More
precisely, we say that a deformation $(i,\phi)$ of $(X,x)$ over
$(S,s)$ is \emph{complete} if any other deformation over some base
space $(T,t)$ can be induced from $(i,\phi)$ by some base change
$\varphi:(T,t) \to (S,s)$. A complete deformation is called
\emph{versal} if for any deformation of $(X,x)$ over some subgerm
$(T',t) \subset (T,t)$ induced by some base change $\varphi ':
(T',t) \to (S,s), \, \varphi$ can be chosen in such a way that it
extends $\varphi '$. We will now give the formal definitions.
\begin{definition} {\rm (1) A deformation $(X,x)\overset{i}{\hookrightarrow}
(\mathscr{X},x)\overset{\phi}{\rightarrow}(S,s)$ of $(X,x)$ is
called \emph{complete} if, for any deformation
$(j,\psi):(X,x)\overset{j}{\hookrightarrow}
(\mathscr{Y},y)\overset{\psi}{\rightarrow}(T,t)$ of $(X,x)$ there
exists a morphism $\varphi:(T,t)\to (S,s)$ such that $(j,\psi)$ is
isomorphic to the induced deformation
$(\varphi^*i,\varphi^*\phi)$.\\
(2) A deformation $(i,\phi)$ is called \emph{versal}  if, for any
deformation $(j,\psi)$ as above
 the following holds: for any closed embedding $k:(T',t)\hookrightarrow
(T,t)$ of complex space germs and any  morphism
$\varphi':(T',t)\to (S,s)$ such that $(k^*j,k^*\psi)$ is
isomorphic to the induced deformation
$(\varphi'^*i,\varphi'^*\phi)$, there exists a
 morphism $\varphi:(T,t)\to (S,s)$
such that\\
(i) $\varphi \circ k = \varphi'$ and \\
(ii)$(j,\psi)$ is isomorphic to the induced deformation
$(\varphi^*i,\varphi^*\phi)$.

This definition can be illustrated by the following commutative
diagram:
$$\quad (X,x)\quad$$
 $$k^*j\swarrow \quad \downarrow j \quad \searrow i $$
$$k^*(\mathscr{Y},y)\quad \hookrightarrow\quad (\mathscr{Y},y)\quad \dashrightarrow (\mathscr{X},x) $$
 $$k^*\psi \downarrow \quad\quad\quad\quad \quad \psi \downarrow \quad\quad\quad \quad \downarrow \phi $$
$$(T',t)\quad \stackrel {k}{\hookrightarrow}\quad\quad(T,t)\quad \stackrel {\varphi}{\dashrightarrow}(S,s). $$

(3) A versal deformation is called \emph{semiuniversal or
miniversal} if the Zariski tangent map $T(\varphi):T_{(T,t)}\to
T_{(S,s)}$ is uniquely defined by $(i,\phi)$ and $(j,\psi)$. }
\end{definition}
Our definition of versality is very restrictive. The deformations
that we call complete are sometimes called versal in the
literature. For example, in \cite{AGV} the authors call our
complete deformation "versal" and our versal deformation
"infinitesimally versal".

Now we introduce a notion which is weaker than completeness, but
still strong enough for many applications.
\begin{definition}{\rm A deformation $(X,x)\overset{i}{\hookrightarrow}
(\mathscr{X},x)\overset{\phi}{\rightarrow}(S,s)$ of $(X,x)$ is
called \emph{1-complete} if, for any one parametric deformation
$(j,\psi):(X,x)\overset{j}{\hookrightarrow}
(\mathscr{Y},y)\overset{\psi}{\rightarrow}(\C,0)$ of $(X,x)$ there
exists a morphism $\varphi:(\C,0)\to (S,s)$ such that $(j,\psi)$
is isomorphic to the induced deformation
$(\varphi^*i,\varphi^*\phi)$}.
\end{definition}

An arbitrary complex space germ may not have a versal deformation.
It is a fundamental theorem of Grauert, that for isolated
singularity a semiuniversal deformation exists.

\begin{theorem}[Grauert, 1972]. Any complex space germ with isolated
singularity has a semiuniversal deformation.
\end{theorem}
\emph{Proof.} See \cite{Gra}.

The following two statements describe the connection between
equisingular families and versal deformations.

\begin{theorem} \label{TjurinaDeform}
Let $(X,0)\subset(\mathbb{C}^n,0)$ be an isolated singularity
defined by $f \in \mathcal{O}_{\C^n,0}$ and $g_1,\ldots,g_{\tau}
\in \mathcal{O}_{\C^n,0}$ be a basis of the Tjurina algebra $T_f$.
If we set $$F(x,t):=f(x)+\sum_{j=1}^{\tau}t_jg_j(x), \quad
(\mathscr{X},0):= V(F) \subset (\C^n \times \C^{\tau},0),$$ then
$(X,0) \hookrightarrow (\mathscr{X},0) \overset{\phi}{\to}
(\C^{\tau},0)$, where $\phi$ is the second projection, is a
semiuniversal deformation of $(X,0)$.
\end{theorem}
\emph{Proof.} See \cite{GLS}, corollary II.1.17.

\begin{corollary} \label{VerDef} Let $(H,z)$ be a germ of a projective hypersurface with one isolated singularity $z$.
Suppose that the equianalytic family of $(H,z)$ is T-smooth at
$(H,z)$. Then the linear system $|H|$ induces a versal deformation
of $(H,z)$.
\end{corollary}

We finish this section with a version of the classical curve
selection lemma that will be used below.
\begin{lemma} \label{CurSel}
Let $X$ be an algebraic variety over $\C$. Let $U \subset X$ be a
Zariski open subset and $x$ be a point in the closure of $U$. Then
there exists a morphism of analytic germs $\phi:(\C,0)\to (X,x)$
such that $\phi(\C\setminus \{0\})\subset U$.
\end{lemma}
\emph{Proof.} A basic lemma (see, for example, \cite{Kem}, lemma
7.2.1) says that there exists a smooth curve $C$ and a morphism
$\nu:C \to \overline{U}$, such that $\nu^{-1}(U)$ is non-empty and
$x$ is contained in the image of $\nu$. Denote
$Z=\nu^{-1}(\overline{U} \setminus U)$ . Then $Z$ is a closed
subset of $C$ and hence consists of a finite number of points. Now
take any point $z \in \nu^{-1}(\{x\}) \subset C$. It has a
neighborhood which does not contain other points of $Z$. Since $C$
is smooth, the analytic germ $(C,z)$ is isomorphic to $(\C,0)$.
Hence $\nu$ defines the required morphism. \proofend

\subsection{Notions of computer
algebra}\label{SecComp}

In our computations we want to use the methods of computer
algebra. Here we introduce the basic notions, that will be widely
used in our proofs. A more detailed description of these notions
can be found also in \cite{GP}.
\subsubsection{Monomial orderings}
\label{MonOrd}
\begin{definition}\label{d7}{\rm A \emph{monomial ordering} is a
total (or linear) ordering $>$ on the set of monomials
$Mon_{n}=\{x^{\alpha}\mid\alpha\in\Z_{\geq 0}^{n}\}$ in $n$
variables satisfying $$x^{\alp}>x^{\bet}\Rightarrow
x^{\gam}x^{\alp}>x^{\gam}x^{\bet}$$  for all
$\alp,\bet,\gam\in\Z_{\geq 0}^{n}$. We say also $>$ is a monomial
ordering on $A[x_{1},...,x_{n}]$, where $A$ is any ring, meaning
that $>$ is a monomial ordering on $Mon_{n}$.}
\end{definition}
We identify $Mon_{n}$ with $\Z_{\geq 0}^{n}$, and then a monomial
ordering is a total ordering on $\Z_{\geq 0}^{n}$, which is
compatible with the semigroup structure on $\Z_{\geq 0}^{n}$ given
by addition. From a practical point of view, a monomial ordering
$>$ allows us to write a polynomial $f\in K[x]$ in a unique
ordered way as
$$f=a_{\alp}x^{\alp}+a_{\bet}x^{\bet}+...+a_{\gam}x^{\gam},$$
with $x^{\alp}>x^{\bet}>...>x^{\gam}$, where no coefficient is
zero.\\
The most important distinction is between global and local
orderings.
\begin{definition}\label{d9}{\rm Let $>$ be a monomial ordering
on $\{x^{\alp}\mid\alp\in\Z_{\geq 0}^{n}\}$. \\
(1) $>$ is called a \emph{global} ordering if $x^{\alp}>1$ for all
$\alp\neq(0,...,0)$, \\
(2) $>$ is called a \emph{local} ordering if $x^{\alp}<1$ for all
$\alp\neq(0,...,0)$}.
\end{definition}
Important examples of monomial orderings are:
\begin{example}\label{mo}{(monomial orderings).
\rm In the following examples we fix an enumeration
$x_{1},...,x_{n}$ of the variables, any other enumeration leads to
a different ordering.\\
(1) \textsc{Global orderings}

(i) \emph{Lexicographical ordering} $>_{lp}$
$$x^{\alp}>_{lp}x^{\bet}:\Leftrightarrow \exists 1\leq i\leq n:
\alp_{1}=\bet_{1},...,\alp_{i-1}=\bet_{i-1},\alp_{i}>\bet_{i},$$

(ii)\emph{ Degree lexicographical ordering} $>_{Dp}$
$$x^{\alp}>_{Dp}x^{\bet}:\Leftrightarrow deg x^{\alp}>deg x^{\bet}$$ $$or\quad (deg x^{\alp}=deg x^{\bet}\quad and
\quad \exists 1\leq i\leq n:
\alp_{1}=\bet_{1},...,\alp_{i-1}=\bet_{i-1},\alp_{i}>\bet_{i}).$$

(iii) \emph{Weighted degree lexicographical ordering}
$Wp(\om_{1},...,\om_{n})$ \\
Given a vector $\omega=(\om_{1},...,\om_{n})$ of integers, we
define the \emph{weighted degree} of $x^{\alp}$ by $
deg_{\om}(x^{\alp}):=\langle \om,\alp \rangle
:=\om_{1}\alp_{1}+\cdots+\om_{n}\alp_{n},$ that is, the variable
$x_{i}$ has degree $\om_{i}$. For a polynomial
$f=\sum_{\alp}a_{\alp}x^{\alp}$, we define the \emph{weighted
degree},
$$deg_{\om}(f):=max\{deg_{\om}(x^{\alp})\mid a_{\alp}\neq 0\}.$$

Using the weighted degree in (ii), with all $\om_{i}>0$, instead
of the usual degree, we obtain the weighted degree
lexicographical ordering, $Wp(\om_{1},...,\om_{n})$.\\\\
(2) \textsc{Local orderings}

(i) \emph{Negative lexicographical ordering} $>_{ls}$
$$x^{\alp}>_{ls}x^{\bet}:\Leftrightarrow \exists 1\leq i\leq n:
\alp_{1}=\bet_{1},...,\alp_{i-1}=\bet_{i-1},\alp_{i}<\bet_{i},$$

(ii) \emph{Negative degree lexicographical ordering} $>_{Ds}$:
$$x^{\alp}>_{Ds}x^{\bet}:\Leftrightarrow deg x^{\alp}<deg x^{\bet}$$ $$or\quad (deg x^{\alp}=deg x^{\bet}\quad and
\quad\exists 1\leq i\leq n:
\alp_{1}=\bet_{1},...,\alp_{i-1}=\bet_{i-1},\alp_{i}>\bet_{i}).$$

(iii) \emph{Negative weighted degree lexicographical ordering}
$Ws(\om_{1},...,\om_{n})$ is a weighted version of the last
ordering.}
\end{example}

\begin{definition}\label{d8}{\rm Let $>$ be a fixed monomial
ordering. Let $f\in K[x], f\neq 0$. Then $f$ can be written in a
unique way as a sum of non-zero terms
$$f=a_{\alp}x^{\alp}+a_{\bet}x^{\bet}+...+a_{\gam}x^{\gam},\quad
x^{\alp}>x^{\bet}>...>x^{\gam},$$ and
$a_{\alp},a_{\bet},...,a_{\gam}\in K.$ We define:\\
(1) $LM(f):=x^{\alp},$ the \emph{leading monomial} of $f$, \\
(2) $LE(f):=\alp$, the \emph{leading exponent} of $f$, \\
(3) $LT(f):=a_{\alp}x^{\alp},$ the \emph{leading term} of $f$, \\
(4) $LC(f):=a_{\alp},$ the \emph{leading coefficient} of $f$, \\
(5) $tail(f):=f-LT(f)=a_{\bet}x^{\bet}+...+a_{\gam}x^{\gam},$ the
tail of $f$.}
\end{definition}

\begin{definition}\label{d10}{\rm For any monomial ordering $>$ on
$Mon(x_{1},...,x_{n})$, we define the ring $K[x]_{>}$ associated
to $K[x]$ and $>$ by
$$K[x]_{>}:=\{\frac{f}{u}\mid f,u \in K[x],\quad
LM(u)=1\}$$ }
\end{definition}
Note that $K[x]_{>}=K[x]$ if and only if $>$ is global and
$K[x]_{>}=K[x]_{\langle x_{1},...,x_{n}\rangle}$ if and only if
$>$ is local.

\subsubsection{Normal form}
 Let $>$ be a monomial ordering and let
 $R=K[x_{1},...,x_{n}]_{>}$ (see Definition \ref{d10} above). For any subset $G\subset R$ define the
ideal
 $$L_{>}(G):=L(G):=\langle LM(g)|g\in G\setminus \{0\}\rangle_{K[x]}.$$
 $L(G)\subset K[x]$ is called the \emph{leading ideal } of G.
 Note that if $I$ is an ideal, then $L(I)$ is the ideal generated
 by all leading monomials of all elements of $I$ and not only by
 the leading monomials of a given set of generators of $I$.

\begin{definition}\label{d11}{\rm Let $\mathcal{G}$ denote the set
of all finite subsets $G\subset R$. A map
$$NF:R\times\mathcal{G}\rightarrow R,\quad(f,G)\mapsto NF(f\mid G),$$
is called a \emph{normal form} on $R$ if, for all $f\in R$ and $G\in \mathcal{G}$,\\
(0) $NF(0|G)=0$;\\
(1) $NF(f|G)\neq 0 \Rightarrow LM(NF(f|G))\notin L(G)$;\\
(2) if $G=\{g_{1},...,g_{s}\}$, then $r:=f-NF(f|G)$ has a standard
representation with respect to $G$, that is, either $r=0$, or
$$r=\sum_{i=1}^{s}a_{i}g_{i},\quad a_{i}\in R,$$
satisfying $LM(f)\geq LM(a_{i}g_{i})$ for all $i$ such that
$a_{i}g_{i}\neq 0$.\\
$NF$ is called a \emph{reduced normal form}, if, moreover,
$NF(f|G)$ is reduced with respect to $G$, i.e. no monomial of the
power series
expansion of $NF(f|G)$ is contained in $L(G)$. \\
As we can see from the definition, $NF(f|G)=0$ if and only if
$f\in\langle G\rangle$.}

\end{definition}

\subsubsection{\textsc{RedNFBuchberger} algorithm for computation of
normal form} \label{RedNFB}

\begin{algorithm}\label{a1}{\rm (\textsc{redNFBuchberger} algorithm)\\
Assume that $>$ is a global monomial ordering.\\
\textsf{Input:} $f\in K[x]$, $G\in \mathcal{G}$\\
\textsf{Output:} $p\in K[x]$, a reduced normal form of $f$ with
respect to $G$.}

 \begin{enumerate}\item {\rm $p:=0$; $h:=f$;
 \item while $(h\neq 0)$
 \begin{enumerate}
 \item while ($h\neq 0$ and $G_{h}:=\{g\in G | LM(g)$ divides
 $LM(h)\}\neq \emptyset)$ \label{InnLoop}
        \\ \quad$\{$ choose any $g\in G_{h}$;
      \\ \quad $h:=h-(LT(h)/LT(g))\cdot g\}$
      \item if $(h\neq 0)$
  \\ \quad $\{p:=p+LT(h)$;
  \\ \quad $h:=tail(h)\}$;
\end{enumerate}
 \item return $p/LC(p)$;}
\end{enumerate}

\end{algorithm}

The algorithm works in the following way: the inner loop
(\ref{InnLoop}) runs until it meets an "obstruction", i.e. the
first monomial that isn't divisible by the leading monomial of any
member of $G$. When the inner loop (\ref{InnLoop}) stops, $h$
stores a normal form of $f$. To make this normal form reduced, we
add the leading term of $h$, i.e. the "obstruction", to $p$ and
continue working with the tail of $h$ in the same way.

Note that any specific choice of "any $g\in G_{h}$" can give a
different normal form function.  For proof of correctness of the
algorithm see \cite{GP}, section 1.6 algorithms 1.6.10 and 1.6.11.

\subsubsection{Highest corner}
\begin{definition}\label{d13}{\rm Let $>$ be a monomial ordering
on $Mon(x_{1},...,x_{n})$ and let $I\subset
K[x_{1},...,x_{n}]_{>}$ be an ideal. A monomial $m\in
Mon(x_{1},...,x_{n})$ is called the \emph{highest corner} of $I$
(with respect to $>$), denoted by
$HC(I)$, if \\
(1) $m\notin L(I)$;\\
(2) $m'\in Mon(x_{1},...,x_{n})$, $m'<m\Rightarrow m'\in L(I)$.}
\end{definition}

\begin{lemma}\label{L2}
{\rm Let $>$ be a monomial ordering on $Mon(x_{1},...,x_{n})$ and
let $I\subset K[x_{1},...,x_{n}]_{>}$ be an ideal. Let $m$ be a
monomial such that $m'<m$ implies $m'\in L(I)$. Let $f\in
K[x_{1},...,x_{n}]$ such that $LM(f)<m$. Then $f\in I$.}
\end{lemma}
\emph{Proof.} See \cite{GP} lemma 1.7.13.

\begin{lemma}\label{L1}{\rm Let $>$ be a weighted degree ordering
on $Mon(x_{1},...,x_{n})$. Moreover, let $f_{1},...,f_{k}$ be a
set of generators of the ideal $I\subset K[x_{1},...,x_{n}]_{>}$
such that $J:=\langle LM(f_{1}),...,LM(f_{k})\rangle$ has a
highest corner $m:=HC(J)$ and $f\in K[x_{1},...,x_{n}]_{>}$. If
$LM(f)<HC(J)$ then $f\in I$.} \end{lemma} \emph{Proof.} See
\cite{GP} lemma 1.7.17.

\subsection{Affine coordinates and the stratum $V^U_d$} \label{VdU}

In this section we enter the notion of stratum $V^U_d$ that will
be used in all the proofs, in order to work in affine coordinates.

\begin{definition}
{\rm Let $S$ be an analytic singularity type of projective
hypersurfaces. Fix homogeneous coordinates on $\PP^n$ and consider
the open subset $$U = \{(t_0: \ldots :t_n) | \, t_0 \neq 0\}
\subset \PP^n.$$ We define $V^U_d(S)$ to be the space of all
hypersurfaces of degree $d$ that have a unique singular point
inside $U$ of singularity type $S$.}
\end{definition}

Note that both $V^U_d(S)$ and $V_d(S)$ are open subsets of the
space of all hypersurfaces of degree $d$ that have at least one
isolated singular point of singularity type $S$. Hence for any
hypersurface $H \in V^U_d(S) \cap V_d(S) $, the germs
$V_{d,H}:=(V_d(S),H)$ and $V^U_{d,H}:=(V^U_d(S),H)$ coincide.
Hence we will formulate statements on $V_{d,H}$ and prove them on
$V^U_{d,H}$.

\begin{remark}
{\rm There exists a natural embedding $\nu:V^U_{d}(S)
\hookrightarrow V^U_{d+1}(S)$ defined in the following way. Let
$F(t_0, \ldots ,t_n)$ be an equation of the hypersurface $H\in
V^U_{d}(S)$. Then we define $\nu(H)$ to be the hypersurface
defined by the equation $t_0F(t_0, \ldots ,t_n)=0$. Note that in
the coordinate system $x_i=\frac{t_i}{t_0}$ on $U$, $H$ and
$\nu(H)$ will be given by the same local equation.

Using the above embedding, the scheme theoretic structure on
$V^U_{d}(S)$ can be computed in the following way. By Theorem
\ref{p3} there exists $N$ such that $V^U_{d+N}(S)$ is a smooth
variety. Then $V^U_{d}(S)$ is equal to the scheme-theoretic
intersection $V^U_{d+N}(S) \cap |{\mathcal O}_{\PP^n}(d)|$.}
\end{remark}

\section{Main results} \label{Main}
\setcounter{lemma}{0}

 We start with a generalization of Example
\ref{Gur1} to higher dimensions:

\begin{theorem} \label{QHomN}
Let $H\subset \mathbb{P}^n$, $n\geq 3$ be the projective
hypersurface given  by the equation $\sum \limits
_{i=1}^{n}t_i^{\alpha_i}t_0^{d-\alpha_i}+\sum \limits
_{i=1}^{n}\lambda_i t_i^d=0$ where $d=\sum \limits
_{i=1}^{n}\alpha_i-(2n+1)$  and the $\lambda_i$ are complex
numbers such that $z=(1,0,...,0)$ is the unique singular point of
$H$. Note that generic $\lambda_i$ satisfy this condition.

Let $V_{d,H}$ be the germ at $H$ of the equianalytic family of
$H$. Then for any $\{\alpha_i\}_{i=1}^n$ such that $d\geq \alpha_i
\geq 2 $ for all $i$, $V_{d,H}$ is
non-T-smooth and $h^1(\kj_{Z^{ea}(H)/\mathbb{P}^{n}}(d))=1$.\\
Furthermore:

(i) If $n=4$ and $d=3$ (i.e.
$\alpha_1=\alpha_2=\alpha_3=\alpha_4=3$), the germ $V_{3,H}(S)$ is
a smooth variety of non-expected codimension (one less than
expected).

(ii) Otherwise, the germ $V_{d,H}(S)$ is a reduced irreducible
non-smooth variety of expected codimension which has a smooth
singular locus. Moreover, the germ $V_{d,H}(S)$ has the sectional
singularity type $A_{1}$.
\end{theorem}
(For proof see Section \ref{ProofQHomN})
\begin{remark}{\rm
It can be shown that the exceptional case of this Theorem can be
generalized in the following
way:\\
Let $H \subset \PP^n$ be the hypersurface given by the local
equation $\sum_{i=1}^n t_i^d=0$, $d\geq \max\{3,7-n\}$. Then the
germ $V_{d,H}$ is an orbit of $PGL_{n+1}$ and hence is a smooth
variety of non-expected codimension.}
\end{remark}

In Example \ref{Gur1} and Theorem \ref{QHomN} we have seen several
examples of equianalytic strata of minimal obstructedness which
are non-reduced, or reducible or have unexpected dimension. In
these examples non-reduced families have smooth reduction, all
components of reducible families are smooth and have expected
codimension, and non-smooth families have smooth singular loci.
The first statement follows from minimal obstructedness. We
conjecture that the other two statements hold for general families
of minimal obstructedness of Newton non-degenerate hypersurface
singularities.

Also, we conjecture that if
$h^1(\kj_{Z^{ea}(H')/\mathbb{P}^{n}}(d))$ is constant along the
equianalytic family of a unisingular projective hypersurface $H$,
then the family has smooth reduction. This can be easily proven
for reduced families (see \cite{Gou2}, Proposition 2.4.1).

The next question that naturally arose was the behavior of the
geometric properties of equianalytic families with respect to the
stabilization of the singularities (see Subsection
\ref{SecStEqui}).

We found out that these phenomena are not stable. Namely, if we
add a new variable $x_{n+1}$ to the space and $x_{n+1}^2$ to the
local equation of the hypersurface, the equianalytic stratum of
the new hypersurface has the same $h^1$ and $\tau$ but is reduced
irreducible of expected codimension. Apparently, the same is true
for any singularity of minimal obstructedness, though sometimes
more variables and their squares need to be added.

More generally, for any hypersurface singularity with $h^1>0$ but
$h^1(2d-2)=0$, the equianalytic stratum obtains an irreducible
component which is reduced of expected dimension after adding
$h^1+1$ squares. The condition $h^1(2d-2)=0$ always holds for
curves. For higher dimensions, $h^1(2d-2)=0$ follows from the
condition $h^1<d-1$.

The following theorem summarizes all that was mentioned above.
\begin{theorem} \label{Squares}
Let $H\subset \PP^n$ be hypersurface of degree $d\geq 3$ with the
unique singular point $z=(1:0:\ldots:0)$.
Let $h^1:=h^1(\kj_{Z^{ea}(H)/\mathbb{P}^{n}}(d))$ and
$\tau:=degZ^{ea}(H)$. Let $F_0$ be the equation of $H$. Suppose
$h^1>0$.

For any $m\geq 1$ define $W^m\subset \PP^{n+m}$ to be the
hypersurface given by equation $F_0+\sum_{j=1}^m
t_{n+j}^2t_0^{d-2}$ and $z_m$ be the point $(1:0:\ldots:0)$.
Denote by $V^{U_m}_{d,W^m}$ the germ at $W^m$ of the family of all
hypersurfaces of degree $d$  that have a unique singular point
inside $U_m = \{(t_0:...:t_{n+m},) |t_0 \neq 0\} \subset
\PP^{n+m}$, and are analytically equivalent to $(W^m,z)$ near the
singular
point. Then: \\\\
(a) $h^1(\kj_{Z^{ea}(W^m)/\mathbb{P}^{n+m}}(d))=h^1$ and $degZ^{ea}(W^m)=\tau$. \\\\
(b) If $h^1(\kj_{Z^{ea}(H)/\mathbb{P}^{n}}(2d-2))=0$ then the
germs $V^{U_m}_{d,W^m}$ for $m\geq h^1+1$ have a reduced component
of
expected dimension.  \\\\
(c) If $H$ is a plane curve then already
$h^1(\kj_{Z^{ea}(H)/\mathbb{P}^{2}}(2d-4))=0$ and hence for $m\geq
h^1+1$ the germs $V^{U_m}_{d,W^m}$ have a reduced component of
expected dimension. \\\\
(d) If $h^1<d-1$ then
$h^1(\kj_{Z^{ea}(H)/\mathbb{P}^{n}}(2d-2))=0$ and hence for $m\geq
h^1+1$ the germs $V^{U_m}_{d,W^m}$ have a reduced component of
expected dimension. \\\\
(e) If $h^1=1$ then the germs $V^{U_m}_{d,W^m}$ are non-smooth of
expected dimension for all $m\geq 1$, reduced for $m \geq
\max\{1,5-d\}$ and irreducible for $m \geq \max\{1,6-d\}$.
\end{theorem}
(For the proof see Section \ref{ProofSquares})

\begin{remark}{\rm
1) Statement (c) is not always true for $n\geq 3$. Consider, for
example, $H$ given by the local equation $\sum x_i^d=0$. Then
$h^1(\kj_{Z^{ea}(H)/\mathbb{P}^{n}}(k))>0$ for $k< n(d-2)$.\\
2) It can be proven that if
$h^1(\kj_{Z^{ea}(H)/\mathbb{P}^{n}}(2d-2))=0$ and
$h^1(\kj_{Z^{ea}(H)/\mathbb{P}^{n}}(d+1))=h^1(\kj_{Z^{ea}(H)/\mathbb{P}^{n}}(d))-1$
then for $m \geq 2$ the germs $V^{U_m}_{d,W^m}$ have a reduced
component of expected dimension. \\
3) If $F_0=\sum_{i=1}^nt_i^{\alp_i}t_0^{d-\alp_i}$ and
$h^1(\kj_{Z^{ea}(H)/\mathbb{P}^{n}}(t))=0$ for some $d<t\leq
2d-2$, then for $m \geq \lceil h^1\frac{(t-d+1)}{2(d-2)}\rceil+1$
the germs $V^{U_m}_{d,W^m}$ have a reduced component of expected
dimension. In particular, if
$F_0=\sum_{i=1}^nt_i^{\alp_i}t_0^{d-\alp_i}$ and
$h^1(\kj_{Z^{ea}(H)/\mathbb{P}^{n}}(2d-3))=0$ then for $m \geq
\lceil \frac{h^1}{2}\rceil+1$ the germs $V^{U_m}_{d,W^m}$ have a
reduced component of expected dimension.  Together with (2) that
implies that if $F_0$ is canonical quasihomogeneous and $h^1<d-1$
then the germs $V_{d,W^m}(S^m)$ have a reduced component of
expected dimension for any  $m \geq \lceil
\frac{h^1}{2}\rceil+1$.\\ For proof see \cite{Gou2}, Theorem
2.2.1(f).}
\end{remark}

\subsection{Deformation theoretic meaning} \label{SecOurDef} In
this subsection we give a deformation theoretic interpretation to
our results. Since the proofs of the statements here are shorter
and less technical, we give them right after the statements.

\begin{theorem}\label{DefMeaning} Let $H\in\PP^n$ be a projective hypersurface of degree
$d$ with the unique singular point $z$. Suppose that the
equianalytic stratum germ $V_{d,H}$ has a reduced component $R$ of
expected dimension. Then the deformation of $H$ induced by the
linear system $|H|$ is 1-complete.
\end{theorem}

The proof of the theorem is based on the following observation: at
every smooth point $H'$ of $R$, the stratum $V_{d,H}$ is
$T$-smooth. Hence the deformation of $H'$ induced by the linear
system $|H|$ is versal, and hence any 1-parametric deformation
$(H,z)\hookrightarrow (\mathscr{X},x)\rightarrow(\C,0)$ of $(H,z)$
can be induced from it by a map $\psi _{H'}: (\C,0) \to (|H|,H')$.
By the curve selection lemma, there exists a map $\phi:(\C,0) \to
(R,H)$ such that all points except $0$ are mapped to non-singular
points. Now we define the requested map $\varphi : (\C,0) \to
(|H|,H)$ by $\varphi(t):= \psi_{\phi(t)}(t)$ for $t \neq 0$ and
$\varphi(0)=H$.

Now we give a precise proof, which includes the description how to
choose the maps $\psi_{\phi(t)}$ analytically.

\emph{Proof.} Denote $U=R\setminus Sing(V_{d})$ where
$Sing(V_{d})$ is the singular locus of $V_{d}$. Let $\tau$ be the
Tjurina number of $(H,z)$. Consider the coincidence variety
$$Z:=\{(H',W)|H' \in R \text{, } W \text{ is a } \tau- \text{
dimensional affine subspace of }|H| \text{ and } H'\in W\}.$$ Let
$Y \subset Z$ be the open subset defined by
$$Y :=\{(H',W)\in Z | H' \in U \text{ and } W \text{ is transversal to } R \text{ at
} H'\}.$$ By the curve selection lemma (Lemma \ref{CurSel}), there
exists a morphism of analytic germs $\phi:(\C,0)\to (Z,(H,W))$
(for some $\tau$-dimensional subspace $W$) such that
$\phi(\C\setminus \{0\})\subset Y$.

Now let $(H,z)\hookrightarrow (\mathscr{X},x)\rightarrow(\C,0)$ be
a one-parametric deformation of $(H,z)$. Let $Ta_H$ denote the
Tjurina algebra of $(H,z)$ and let $(H,z)\hookrightarrow
(\mathscr{Y},y)\rightarrow(Ta_H,0)$ be the semiuniversal
deformation over it described in Theorem \ref{TjurinaDeform}.
Since this deformation is semiuniversal, there exists a morphism
$\psi:(\C,0) \to (Ta_H,0)$ such that $\psi^*(\mathscr{Y},y)=
(\mathscr{X},x)$. Note that the monomial basis of the Tjurina
algebra of $(H,z)$ is also a basis of Tjurina algebras in a
neighborhood of $H$. Thus we identify those Tjurina algebras as
vector spaces.

For any point $(H',W) \in Y$, the factor morphism $T_{H'}|H| \to
Ta_{H'} \cong Ta_H$ defines an isomorphism $p_{H',W}:W \cong
Ta_H$, since $W$ is transversal to the kernel of the factor
morphism, which is the tangent space to $R$ at $H'$.
This defines a morphism $\Psi: (\C,0) \times Y \to |H|$ by
$$\Psi(t,(H',W)):= p_{H',W}^{-1}(\psi(t)). $$

Now, we define $\varphi :(\C,0) \setminus 0 \to (|H|,H)$ by
$\varphi:= \Psi \circ (Id \times \phi)$, and extend it to 0 by
$\varphi(0):=H$. \proofend
%
\begin{corollary}
Let $H\in \PP^n$ be a unisingular hypersurface of degree $d\geq 3$
defined by the equation $\sum \limits
_{i=1}^{n}t_i^{\alpha_i}t_0^{d-\alpha_i}+\sum \limits
_{i=1}^{n}\lambda_i t_i^d=0$. Suppose that $H$ has one isolated
singularity and $d+1=\sum \limits _{i=1}^{n}(\alpha_i-2)$. Then,
unless $n=2$, $d\leq 6$ or $n=4$, $d=3$, the deformation of $H$
induced by the linear system $|H|$ is 1-complete.
\end{corollary}

\begin{corollary} \label{SquareCor}
Let $H\subset \PP^n$ be a hypersurface of degree $d\geq 3$ with
the unique singular point $z=(1:0:\ldots:0)$. Let $F_0$ be the
equation of $H$. For any $m\geq 1$ define $W^m\subset \PP^{n+m}$
to be the hypersurface given by equation $F_0+\sum_{j=1}^m
t_{n+j}^2t_0^{d-2}$ and $z_m$ be the point $(1:0:\ldots:0)$.

Suppose that $h^1(\kj_{Z^{ea}(H)/\mathbb{P}^{n}}(2d-2))=0$. Then
for $m \geq h^1+1$ the deformation of $W^m$ induced by the linear
system $|W^m|$ is 1-complete.
\end{corollary}

Let us now demonstrate one known application of 1-completeness.
Suppose that we want to construct a hypersurface of degree $d$
having $m$ isolated singular points of prescribed analytic
singularity types $S_1$,...,$S_m$. Suppose that we can construct a
hypersurface $H$ with unique more complicated singularity that
splits to singularities of the given types $S_i$ after a
one-parameter deformation by hypersurfaces of higher degrees. If
the deformation of $H$ induced by the linear system $|H|$ is
1-complete, there exists a deformation of $H$ by hypersurfaces
from $|H|$ which contains the desired hypersurfaces.

\begin{proposition}
Let $H \subset \PP^n$ be a hypersurface of degree $d$ with one
isolated singular point $z$ of the analytic singularity type $S$.
Let $(H,z)\overset{i}{\hookrightarrow}
(\mathscr{X},x)\overset{\phi}{\rightarrow}(|H|,H)$ be the
deformation of $(H,z)$ induced by the linear system $|H|$. Suppose
that it is 1-complete. Let $S_1$,..,$S_m$ be analytic singularity
types. Suppose also that there exists a one-parameter deformation
$(j,\psi):(H,z)\overset{j}{\hookrightarrow}
(\mathscr{Y},y)\overset{\psi}{\rightarrow}(\mathbb{C},0)$ of
$(H,z)$ that includes hypersurfaces having $m$ singularities of
types $S_1$,...,$S_m$. Then there exists a one-parameter
deformation of $(H,z)$ consisting of hypersurfaces of degree $d$
that includes hypersurfaces having $m$ singularities of types
$S_1$,...,$S_m$.
\end{proposition}
\emph{Proof.} Since the deformation of $H$ induced by $|H|$ is
1-complete, there exists a morphism $\varphi:(\mathbb{C},0)\to
(|H|,H)$ such that $(j,\psi)$ is isomorphic to the induced
deformation $(\varphi^*i,\varphi^*\phi)$. Hence the induced
deformation $(\varphi^*i,\varphi^*\phi)$ includes hypersurfaces
having $m$ singularities of types $S_1$,...,$S_m$. On the other
hand, the deformation $(i,\phi)$ consists of hypersurfaces of
degree $d$, and hence the induced deformation
$(\varphi^*i,\varphi^*\phi)$ also consists of hypersurfaces of
degree $d$. \proofend

\section{Proof of the theorem on quasihomogeneous
hypersurface singularities}\label{ProofQHomN}
\setcounter{lemma}{0}
This section is dedicated to the proof of Theorem \ref{QHomN}.
%
%
%

\subsection{The structure of the proof}\label{Struct2}
First of all we pass to affine coordinates $x_i=\frac{t_i}{t_0}$.
In these coordinates, $H$ is given by the local equation
$f=\sum_{i=1}^n x_i^{\alpha_i}+\sum_{i=1}^n \lambda_i x_i^d=0$. It
is a semiquasihomogeneous polynomial with non-degenerate
quasihomogeneous part $g=\sum_{i=1}^n (1+\delta_{\alp_i,d}
\lambda_i) x_i^{\alpha_i}$. It is easy to see that $f\in j(f)$ and
$j(f)=j(g)$. Hence $f$ and $g$ have the same Tjurina ideal and
Tjurina algebra, $f\overset{c}{\sim}g$ and $T_f=M_f$.

Then we prove the equality
$h^1(\kj_{Z^{ea}(H,z)/\mathbb{P}^{n}}(d))=1$ (see Subsection
\ref{Hyperh1=1}). Next (in Subsection \ref{SubNot}) we switch to
substratum germ $V_{d,H}^{0,0}$ of $V_{d,H}$ consisting of
hypersurfaces given by polynomials of the form
\begin{multline} F(x_1,...,x_n)=\sum_{i=1}^{n}x_i^{\alp_i}+\sum_{i=1}^{n}\lambda_i x_i^d+\sum_{I\in \mathcal{D}}
 a_{I}x^{I},\\ \text{ where }\mathcal{D}=(\{I\in \Z_{\geq 0}^n | \alp_n \leq |I|\leq d\}\setminus
 \bigcup_{1\leq i\neq j \leq
 n}\{(0,...,0,\alp_i-1,0,...,0,\underset{j}{1},0,...,0)\})
 \setminus\\
\setminus\bigcup_{1\leq i \leq n}\{(0,...,0,\alp_i,0,...,0)\}.
\nonumber\end{multline}

We prove that this substratum is transversal to the orbits of the
group of affine transformations of $\C^n$.

First, we consider the case $\alp_1<2\alp_n$ (Subsection
\ref{HyperCase1}). In this case for any hypersurface $H$ which
lies in the stratum germ there exists an affine coordinate change
s.t. the equation of $H$ in the new coordinates does not include
any terms that lie below the Newton polytope $\Delta(f)$, and has the same terms laying on $\Delta(f)$ as $f$.

Let $F=f+f_1$ where $f_1$ is a polynomial which has no terms below
and on $\Delta(f)$. We claim that $F$ is contact equivalent to $f$
if and only if $F\in j(F)=\langle F_{x_1},\dots,F_{x_n}\rangle$.

One direction is obvious: if they are equivalent then they have
the same Milnor and Tjurina numbers and hence $\mu(F)=\tau(F)$,
i.e. $F \in j(F)$. To prove the other direction we use Saito
theorem (Theorem \ref{Saito}). It says that if $F \in j(F)$ then
there exists a quasihomogeneous polynomial $h$ and a coordinate
change $\phi$ that maps $h$ to $F$. Then the linear part of $\phi$
will map $h$ to the quasihomogeneous part of $F$, which is $g$.
Therefore, $h$ and $g$ are contact equivalent and hence $F$ and
$f$ are contact equivalent.
So the hypersurface $H_F$ belongs to $V_{d,H}^{0,0}$ if and only
if $F\in \langle F_{x_1},\dots,F_{x_n}\rangle$.

We check that condition using a computer algebra algorithm
(Algorithm \ref{OurRedNFBuch}). In this way we obtain a system of
equations on $V_{d,H}^{0,0}$. In case $d=3,\, n=4$ the substratum
germ consists of one point. We show that otherwise the system
consists of a subsystem having a diagonal linear part, and one
more equation with quadratic principle part of rank $\geq 3$.

In the case of $\alp_1\geq 2\alp_n$, there are hypersurfaces $H_F$
in $V_{d,H}^{0,0}$ whose equations include some terms below or on
the Newton polytope $\Delta(f)$. For every such polynomial $F$, we
pass to new coordinates in which $F$ has no terms below and on
$\Delta(f)$, write equations on the coefficients of $F$ in the new
coordinates and express new coefficients through the old ones.

Again we check that the obtained system consists of a subsystem
having diagonal linear part, and one more equation with quadratic
principle part of rank $\geq 3$. This is done in Subsection
\ref{HyperCase2}.

We show that in both cases the last equation lies in the ideal
generated by elements that appear in its quadratic part. We deduce
from this fact the smoothness of the singular locus.
\subsection{Proof that
$h^1(\kj_{Z^{ea}(H,z)/\mathbb{P}^{n}}(d))=1$} \label{Hyperh1=1}
Suppose, for convenience, $\alp_1\geq \alp_2\geq...\geq \alp_n
\geq 2$.

First, we pass to affine coordinates $x_i=\frac{t_i}{t_0}$. In
these coordinates, $H$ is given by local equation $f=\sum
x_i^{\alpha_i}+\sum \lambda_i x_i^d$. It is a semiquasihomogeneous
polynomial with non-degenerate quasihomogeneous part
$g=\sum_{i=1}^n (1+\delta_{\alp_i,d} \lambda_i) x_i^{\alpha_i}$.
It is easy to see that $f \in j(f)$ and $j(f)=j(g)$. Hence $f$ and
$g$ have the same Tjurina ideal and Tjurina algebra, which also
coincides with their Milnor algebras. By Mather-Yau theorem this
implies that $f$ and $g$ are contact equivalent and hence belong
to the same stratum.

The polynomial $g$ is quasihomogeneous of type
$(1/\alp_1,\dots,1/\alp_n;1)$ hence $(H,z)$ is a quasihomogeneous
hypersurface singularity. The Newton polytope of $g$ is
$$\Delta(g) = \{I\in \Z_{\geq 0}^n|
\sum_{j=1}^n\frac{I_j}{\alp_j}=1\}.$$
\begin{figure}[h]
\setlength{\unitlength}{1cm}
\begin{picture}(18,10)(0,0)
\thicklines\put(4,5){\vector(0,1){5}} \put(4,5){\vector(1,0){5}}
 \put(4,5){\vector(-1,-1){2.5}}
\put(4,9.25){\line(1,-1){4.27}} \put(4,9.25){\line(-1,-4){1.41}}
 \put(2.53,3.50){\line(4,1){5.75}}
\thinlines
\dashline{0.2}(3.25,4.25)(3.25,5.5)
\dashline{0.2}(4.75,4.25)(4.75,5.5) 
\dashline{0.2}(5.5,5)(5.5,6.25) 
\dashline{0.2}(3.25,4.25)(4.75,4.25) 
\dashline{0.2}(4.75,4.25)(5.5,5)
\dashline{0.2}(3.25,5.5)(4,6.25)
\dashline{0.2}(4,6.25)(5.5,6.25)
\dashline{0.2}(5.5,6.25)(4.75,5.5)
\dashline{0.2}(4.75,5.5)(3.25,5.5)

\put(2.5,3.1){$d$} \put(3.23,3.95){$\alp_1-2$}

\put(4.25,5.4){$\mathcal{P}$} \put(4.25,7.4){$T_d$}

\put(8.4,5.1){$d$} \put(5.7,5.1){$\alp_2-2$}

\put(4.15,9.3){$d$} \put(4.15,6.4){$\alp_3-2$}

\thicklines\put(11,5){\vector(0,1){5}} \put(11,5){\vector(1,0){5}}
\put(11,5){\vector(-1,-1){2.5}} \put(11,9.25){\line(1,-1){4.27}}
\put(11,9.25){\line(-1,-4){1.41}}
 \put(9.53,3.50){\line(4,1){5.75}}
\thinlines
\dashline{0.2}(10.25,4.25)(10.25,5.5)
\dashline{0.2}(11.75,4.25)(11.75,5.5) 
\dashline{0.2}(12.5,5)(12.5,6.25) 
\dashline{0.2}(10.25,4.25)(11.75,4.25) 
\dashline{0.2}(11.75,4.25)(12.5,5)
\dashline{0.2}(10.25,5.5)(11,6.25)
\dashline{0.2}(11,6.25)(12.5,6.25)
\dashline{0.2}(12.5,6.25)(11.75,5.5)
\dashline{0.2}(11.75,5.5)(10.25,5.5)

\dashline{0.2}(11,6.75)(13,5) 
\dashline{0.2}(10,4)(13,5) 
\dashline{0.2}(10,4)(11,6.75) 

\put(9.5,3.1){$d$} \put(10,3.8){$\alp_1$}

\put(15.4,5.1){$d$} \put(13.1,5.1){$\alp_2$}

\put(11.15,9.3){$d$} \put(11.15,6.8){$\alp_3$}

\end{picture}
\caption {Newton polytope}\label{NewtonPolytope}
\end{figure}
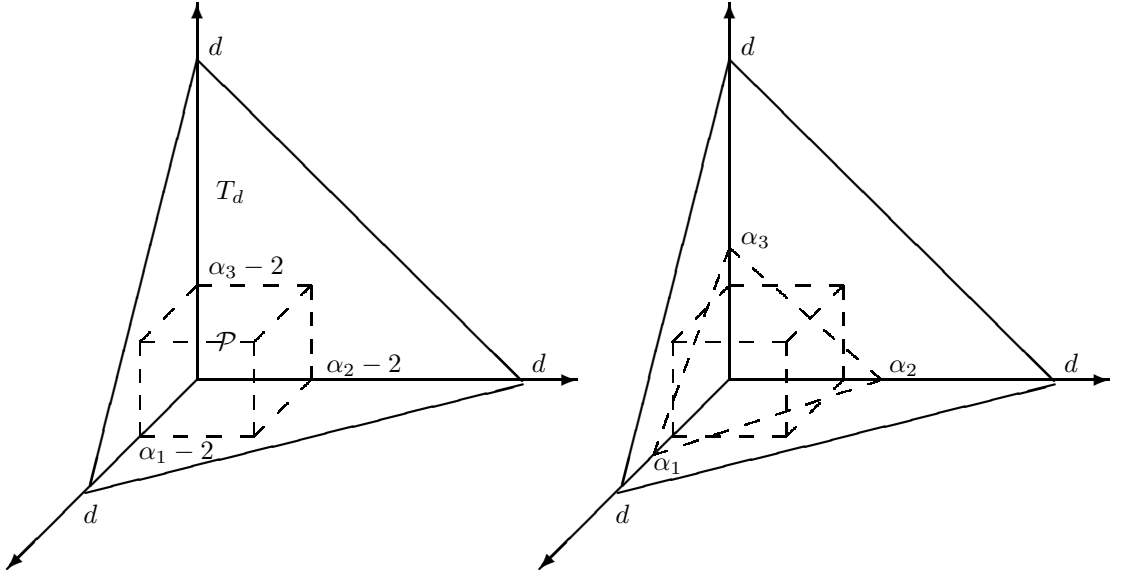

We will now show that $V_{d,H}$ is non-T-smooth at $H$ and
$h^1(\kj_{Z^{ea}(H,z)/\mathbb{P}^{n}}(d))=1$.

The Tjurina algebra of $f$ has a basis $\{x^{I}, I \in
\mathcal{P}\}$, where $\mathcal{P}$ is the parallelepiped
$\mathcal{P}=\{I\in \Z_{\geq 0}^n|\, I_j\leq \alp_j-2$ for all
$0\leq j\leq n\}$. So
$\tau(H,z)=|\mathcal{P}|=\prod_{i=1}^n(\alp_i-1)$, where by
$|\mathcal{P}|$ we denote
the number of integer points in $\mathcal{P}$. \\
Hence $h^0(\mathcal{O}_{Z^{ea}(H,z)})=\tau(H,z)=|\mathcal{P}|$,
$$H^0(\kj_{Z^{ea}(H,z)/\mathbb{P}^{n}}(d))=\{\sum_{I\in
T_d\setminus \mathcal{P}}a_{I}x^I\}$$ where $T_d$ is the simplex
$\{I\in \Z_{\geq 0}^n|\,|I|\leq d \}$ (see Figure
\ref{NewtonPolytope}). That means that
$$h^0(\kj_{Z^{ea}(H,z)/\mathbb{P}^{n}}(d))=|T_d|-|T_d\cap \mathcal{P}|$$
From the exact sequence $$0\rightarrow
H^0(\kj_{Z^{ea}(H,z)/\mathbb{P}^{n}}(d))\rightarrow
H^0(\mathcal{O}_{P^n}(d))\rightarrow
H^0(\mathcal{O}_{Z^{ea}(H,z)})\rightarrow
H^1(\kj_{Z^{ea}(H,z)/\mathbb{P}^{n}}(d))\rightarrow 0$$ we
conclude that

$$h^1(\kj_{Z^{ea}(H,z)/\mathbb{P}^{n}}(d))=h^0(\kj_{Z^{ea}(H,z)/\mathbb{P}^{n}}(d))-
h^0(\mathcal{O}_{P^n}(d))+h^0(\mathcal{O}_{Z^{ea}(H,z)})=$$
$$|T_d|-|T_d\cap \mathcal{P}|-|T_d|+|\mathcal{P}|=|\mathcal{P}\setminus T_d|=1.$$
The same argument shows that
$h^1(\kj_{Z^{ea}(H,z)/\mathbb{P}^{n}}(d+1))=|\mathcal{P}\setminus
T_{d+1}|=0$.
 Thus by Theorem \ref{p3} the germ $V_{d+1,H}$ is T-smooth and the germ $V_{d,H}$ is
non-T-smooth.

Because of minimal obstructedness ($h^1=1$), $V_{d,H}$ may be
either non-smooth of expected codimension or smooth of
non-expected codimension.

Now we would like to find out when it is non-smooth and when it
has non-expected codimension. First we will pass to a more
convenient substratum, which has the same geometric properties.

\subsection{Switch to substratum and notations} \label{SubNot}

First of all let us shift the singularity to the origin. Let $S$
be the singularity type of $(H,z)$. The family $V_{d}^U(S)$ is
invariant under the action of affine transformations of $\C^n$.
Consider the subgroup generated by translations. We switch to the
section $V'_{d}(S)$ of $V_{d}^U(S)$ transversal to orbits of this
group and given by the conditions that the singularity is in the
origin.

In the same way, using the subgroup $GL_n$ consisting of the
linear coordinate changes, we want to reduce to the substratum
$V_{d}^{0,0}(S)$ of $V_{d}'(S)$ consisting of all hypersurfaces
$H_F$ given by polynomials $F$ which also do not include the
monomials $x_i^{\alp_i-1}x_j$ for $i\neq j$, and include the
monomials $x_i^{\alp_i}$ with coefficient 1 if $\alp_i \neq d$ or
with coefficient $1+\lambda_i$ if $\alp_i=d$ . For this purpose we
will prove the following lemma.
\begin{lemma} \label{sub}
$ $\\
(i)  $T_{H}V_{d}'(S)=T_{H}V_{d}^{0,0}(S)\oplus T_{H}GL_nH$\\
(ii) $GL_nV_{d}^{0,0}(S)=V_{d}'(S)$ in a neighborhood of $H$ and
$GL_nH \cap V_{d}^{0,0}(S) = \{H\}$ in a neighborhood of
$H$.\\
The same is true for $V_{d+1}'(S)$:\\
(iii)  $T_{H}V_{d+1}'(S)=T_{H}V_{d+1}^{0,0}(S)\oplus T_{H}GL_nH$\\
(iv) $GL_nV_{d+1}^{0,0}(S)=V_{d+1}'(S)$ in a neighborhood of $H$
and $GL_nH \cap V_{d+1}^{0,0}(S) = \{H\}$ in a neighborhood of
$H$.\\
\end{lemma}
\emph{Proof.} \\
(i): For any point $H_F\in V_{d}'(S)$ and any $0\leq i,j\leq n$
denote by $c_{ij}$ the coefficient of the monomial
$x_i^{\alp_i-1}x_j$  in the polynomial $F-f$. The tangent space to
$V_{d}^{0,0}(S)$ at $H$ is given inside $T_{H}V_{d}'(S)$ by the
equations $c_{ij}=0$. On the other hand $T_{H}GL_nH=Span\{\alp_i
x_i^{\alp_i-1}x_j +\lambda_i d x_i^{d-1}x_j \}$. Hence
$T_{H}V_{d}'(S)=T_{H}V_{d}^{0,0}(S)\oplus
T_{H}GL_nH$. \\
(iii) is proven in the same way.\\
(iv) follows from (iii) since $V_{d+1}^{0,0}(S)$ is smooth (see Section \ref{Hyperh1=1}). \\
(ii) follows from (iv) since $GL_n$ preserves degree:\\
$$GL_nV_{d}^{0,0}(S)=GL_n(V_{d+1}^{0,0}(S)\cap
|\mathcal{O}_{\PP^n}(d)|)=(GL_nV_{d+1}^{0,0}(S))\cap
|\mathcal{O}_{\PP^n}(d)|=$$ $$=V_{d+1}'(S)\cap
|\mathcal{O}_{\PP^n}(d)|=V_{d}'(S).$$ Here,
$|\mathcal{O}_{\PP^n}(d)|=|H|$ is the linear system of
hypersurfaces of degree $d$. Also $GL_nH \cap V_{d}^{0,0}(S)
\subset GL_nH \cap V_{d+1}^{0,0}(S) = \{H\}$ in a neighborhood of
$H$. \proofend \\
From this lemma we see that it is enough to prove our statement
for the germ $V_{d,H}^{0,0}$ of $V_{d}^{0,0}(S)$ at $H$.

Consider now arbitrary hypersurface $H_F\in V_{d,H}^{0,0}$ given
by a polynomial equation $F=0$. Since $H_F$ is obtained from $H$
by a local analytic diffeomorphism and both have their only
singularity
at the origin, $F$ has no terms of degree less than $\alp_n$. 

So we will work with substratum germ $V_{d,H}^{0,0}$ of $V_{d,H}$
consisting of hypersurfaces  given by polynomials of the form 
\begin{multline}\label{polinomF} F(x_1,...,x_n)=\sum_{i=1}^{n}x_i^{\alp_i}+\sum_{i=1}^{n} \lambda_i x_i^d+\sum_{I\in \mathcal{D}}
 a_{I}x^{I}, \\ \text{ where }\mathcal{D}=(\{I\in \Z_{\geq 0}^n | \alp_n \leq |I|\leq d\}\setminus
 \bigcup_{1\leq i\neq j \leq
 n}\{(0,...,0,\alp_i-1,0,...,0,\underset{j}{1},0,...,0)\}) \setminus \\
\setminus\bigcup_{1\leq i \leq n}\{(0,...,0,\alp_i,0,...,0)\}.
\end{multline}
 For convenience, we introduce the following
notations:
\\a) let $b_{I}$ be the coefficients of basis monomials above the Newton polytope $\Delta(f)$
, i.e. $b_{I}:=a_{I}$ for $I \in \mathcal{D}$ such that $I_j\leq \alp_j-2$ for all $j$ and $w(I)>1$;\\
b)  let $e_{I}$ be the coefficients of basis monomials below the
Newton polytope $\Delta(f)$ but of degree at least $\alp_n$   
, i.e. $e_{I}:=a_{I}$ for $I \in \mathcal{D}$ such that $I_j\leq \alp_j-2$ for all $j$ and $w(I)\leq 1$;\\
c) let $g_I:=a_{I}$ for $I \in E$ where $E=\{I \in \mathcal{D} |$
$I_j=\alp_j-1$ for some $j$, $I_k\leq \alp_k-2$
for all $k\neq j$ and $I_k>0$ for some $k \neq j\}$;\\
d) let $u_{I}:=a_{I}$ for $I \in \mathcal{D}$ such that ($I_j\geq
\alp_j$ for some $j$) or ($I_j = \alp_j-1$ and $I_k = \alp_k-1$
for some $k \neq
j$);\\
e) let $q_{I}:=a_{I}$ for $I=(0,...,0,\alp_j-1,0,...,0)$ for some
$j$; \\
f) for $I=(i_1,...,,i_{k-1},\alp_k-1,i_{k+1},...i_n)\in E$ denote
$$dual(I):=(\alp_1-2-i_1,...,,\alp_{k-1}-2-i_{k-1},\alp_k-1,\alp_{k+1}-2-i_{k+1},...,\alp_{n}-2-i_n).$$
Note that $dual(I)$ also lies in $E$ and $dual(dual(I))=I$.

Let $A=\C[a_I]$ be the algebra of polynomials generated by $a_I$,
$I\in \mathcal{D}$.  Let $m=\langle a_I \rangle$ be the maximal
ideal in $A$ generated by all $a_I$, $G=\langle g_I \rangle$ be
the ideal in $A$ generated by all $g_I$ and $B=\langle b_I
\rangle$ be the ideal in $A$ generated by all $b_I$.

\subsection{Proof of the theorem for the case $\alp_1<2\alp_n$}
\label{HyperCase1}

In this case the Newton polytope $\Delta(f)$ lies below the
hyperplane $|I|=2\alp_n$.

We want to find out for which $\{a_{I}\}$ $H_{F}$ lies in
$V_{d,H}^{0,0}$. Our $F$ doesn't include monomials $x^I$ for $I$
below and on the Newton polytope and satisfying $I_j= \alp_j-1$
for some $j$. Hence, by Corollary \ref{Var}, in order to belong to
our substratum, $F$ should include no terms below and on the
Newton polytope except of $x_i^{\alp_i}$.

Let $F=f+f_1$ where $f_1$ is a polynomial which has no terms below
and on $\Delta(f)$. We claim that $F$ is contact equivalent to $f$
if and only if $F\in j(F)=\langle F_{x_1},\dots,F_{x_n}\rangle$.

One direction is obvious: if they are equivalent then they have
the same Milnor and Tjurina numbers and hence $\mu(F)=\tau(F)$,
i.e. $F \in j(F)$. To prove the other direction we use Saito
theorem (Theorem \ref{Saito}). It says that if $F \in j(F)$ then
there exists a quasihomogeneous polynomial $h$ and a coordinate
change $\phi$ that maps $h$ to $F$. Then the linear part of $\phi$
will map $h$ to the quasihomogeneous part of $F$, which is $g$.
Therefore, $h$ and $g$ are contact equivalent and hence $F$ and
$f$ are contact equivalent.
So the hypersurface $H_F$ belongs to $V_{d,H}^{0,0}$ if and only
if $F\in \langle F_{x_1},\dots,F_{x_n}\rangle$.

In order to check whether $F(x_1,...,x_n)\in \langle
F_{x_1},...,F_{x_n} \rangle$ we use the \textsc{redNFBuchberger}
algorithm (Algorithm \ref{a1}). We refer to a neighborhood of the
origin, hence we consider $F(x_1,...,x_n)$ and $\langle
F_{x_1},...,F_{x_n} \rangle$ in the local ring
$R=\mathbb{C}[x_1,...,x_n]_{\langle x_1,...,x_n\rangle}$. To
compute in this ring, we define a local monomial ordering on
$\mathbb{C}[x_1,...,x_n]$ such that the ring associated to
$\mathbb{C}[x_1,...,x_n]$ and this ordering will be
$\mathbb{C}[x_1,...,x_n]_{\langle x_1,...,x_n\rangle}$. We choose
the negative weighted degree lexicographical ordering with
$w=(1/\alp_1,...,1/\alp_n)$ (see Example \ref{mo}, ordering
(2)(iii)).

 In general, the \textsc{redNFBuchberger} algorithm does
not stop for local orderings. However, in our case we can stop it
manually when the leading monomial of the tail is less than
$x_1^{\alp_1-2}\cdot...\cdot x_n^{\alp_n-2}$. We are allowed to do
that by Lemma \ref{L1}, for $x_1^{\alp_1-2}\cdot...\cdot
x_n^{\alp_n-2}$ is the highest corner of $\langle
LM(F_{x_1}),...,LM(F_{x_n})\rangle=\langle
x_1^{\alp_1-1},...,x_n^{\alp_n-1}\rangle$. Indeed, any monomial
smaller than $x_1^{\alp_1-2}\cdot...\cdot x_n^{\alp_n-2}$ has
degree of $x_j$ bigger than or equal to $\alp_j-1$ for some $j$
and hence lies in $\langle LM(F_{x_1}),...,LM(F_{x_n})\rangle$ and
$x_1^{\alp_1-2}\cdot...\cdot x_n^{\alp_n-2}\notin\langle
LM(F_{x_1}),...,LM(F_{x_n})\rangle$.

There is another explanation why we can stop the algorithm at this
point. Consider $V_{d+1,H}$. It is smooth of expected codimension
(see Subsection \ref{Hyperh1=1}). Therefore $V_{d+1,H}^{0,0}$ is
also smooth and has expected codimension which is equal to the
number of basis elements which lie above the Newton polytope.
Since we have exactly this number of independent equations on this
stage, there will be no more equations. Also since
$V_{d+1,H}^{0,0}$ is smooth, all the equations on it will have
independent linear parts. When we return to $V_{d,H}^{0,0}$, the
linear part of only one of them may vanish.

So  we  rewrite the algorithm in the following way:
\begin{algorithm} \label{OurRedNFBuch} {\rm(Modified  \textsc{redNFBuchberger}
algorithm.)}
\begin{enumerate}\item {\rm $p:=0$, $h:=F$;
 \item while ($h\neq 0$ and $LM(h)\geq x_1^{\alp_1-2}\cdot\ldots\cdot x_n^{\alp_n-2}$)
 \begin{enumerate}
 \item while ($h\neq 0$ and $LM(h)\geq x_1^{\alp_1-2}\cdot\ldots\cdot x_n^{\alp_n-2}$\\
 and exists $i$ such that $LM(F_{x_i})$ divides  $LM(h)$)\\
  $\{h:=h-(LT(h)/LT(F_{x_i}))\cdot F_{x_i}\}$\\
  \item if $(h\neq 0$ and $LM(h)\geq x_1^{\alp_1-2}\cdot\ldots\cdot x_n^{\alp_n-2})$\\
  $\{p:=p+LT(h)$;\\
 $h=tail(h)\}$;
 \end{enumerate}
 \item return $p$;}
\end{enumerate}
\end{algorithm}

As a result, we obtain the normal form
$$NF(F|\langle F_{x_1},...,F_{x_n}
\rangle)=\sum R_{I}(a_{J})x^I,$$ where $x^I$, $I_k\leq \alp_k-2$,
for all $k$ are elements of the basis of algebra
$\mathbb{C}[x_1,...,x_n]/{\langle f_{x_1},...,f_{x_n} \rangle}$
which lie above $\Delta(f)$ and
 $R_{I}(a_{J})$ are polynomials in
$a_{J}$. Hence $F$ belongs to the ideal $\langle
F_{x_1},...,F_{x_n} \rangle$ if and only if all coefficients
$R_{I}(a_{J})=0$.

Thus we obtain a system of equations on $a_{J}$:
\begin{equation} \label{RIzero}
R_{I}(a_{J})=0
\end{equation}
Let us now analyze the $R_I$.
\begin{lemma} \label{RI}
Denote $\psi_{I}(g_{J},u_{K},b_L):= R_I - (1-w(I))b_I$ for $|I|
\leq d$ and $\psi_{I}(g_{J},u_{K},b_L) = -R_I$ for $|I|
=d+1$, i.e. $I = (\alpha_1-2,...,\alpha_n-2)$. Then\\
(i) All $b_L$ that appear in $\psi_{I}$ satisfy $w(L)<w(I)$.\\
(ii) All $\psi_{I}$ are polynomials from $G^2+Bm$. Recall that
$G=\langle g_{I}\rangle$ and $B=\langle b_{I}\rangle$.\\
(iii)$\psi_{(\alpha_1-2,...,\alpha_n-2)}-\sum_{J \in E} A_k\cdot
g_{J}\cdot
 g_{dual(J)}\in Bm+m^3$ where $A_k$ are positive rational numbers.
\end{lemma}
\emph{Proof.} Let us trace the changes of the coefficients of $h$
and $p$ during the algorithm. Denote the coefficient of $x^I$ in
$h$ by $c_I$. In the first step we eliminate the monomials
$x_i^{\alp_i}$ for all $i$. As a result we obtain
$h:=F-\sum_{i=1}^n(1/\alp_i)\cdot x_i\cdot F_{x_i}.$ After this
step $c_I=(1-w(I))a_I$. Note that $c_I$ is non-zero iff $a_I$ is
non-zero, since $w(I)>1$ unless $a_I=0$.

Now let $S$ be the coefficient of the leading monomial of $h$. If it
is a basic monomial, we add the leading term of $h$ to $p$ and
subtract it from $h$. Otherwise, $S_m \geq \alp_m-1$ for some $m$ and
hence we can eliminate this term using $F_{x_m}$. After such
elimination step the change of $h$ is expressed by the formula
$$(h)^{new}
= (h)^{old} - \frac{1-w(S)}{\alp_m} a_Ix_1^{S_1} \cdot ... \cdot
x_m^{S_m-\alp_m+1} \cdot ... \cdot x_n^{S_n}F_{x_m}.$$

Hence the $c_I$ after this step is
$$c_I^{new}= (1-w(I))a_I-\frac{I_m-S_m+\alp_m}{\alp_m}(1-w(S))a_S\cdot
a_{I-S+J^m},$$ where $J^m:= (0,...0,\alp_m,0,..,0)$. If some
coordinate of $I-S+J^m$ is negative, then $a_{I-S+J^m}=0$ and
hence $c_I$ did not change.

Suppose that $I$ is a basic index.  Let us show that after this
step $c_I - (1-w(I))b_I \in G^2+Bm$. We know that $S_m \geq
\alp_m-1$. If $S_k > I_k$ for some $k \neq m$ then $(I-S+J^m)_k<0$
and hence $c_I$ did not change. Hence we can assume $S_k \leq I_k$
for $k \neq m$. Let us consider several cases.
\begin{enumerate}
\item If $S_m \geq \alp_m$ then
$a_{I-S+J^m} \in B$ and hence $a_S\cdot a_{I-S+J^m} \in mB$.

\item If $S_m = \alp_m-1$ and $I_m < \alp_m - 2$ then
 $a_{I-S+J^m} \in B$ and hence
$a_S\cdot a_{I-S+J^m} \in mB$.

\item If $S_m = \alp_m-1$ and $I_m = \alp_m - 2$ then $a_S \in G$
and $a_{I-S+J^m} \in G$ and hence $a_S\cdot a_{I-S+J^m} \in G^2$.
\end{enumerate}

Let us now consider any elimination step of the algorithm at which we eliminate some
$x^S$ using $F_{x_m}$ for some $m$. This is only possible if $S_m
\geq \alpha_m-1$. The change of $h$ in this step is expressed by
the formula $$(h)^{new} = (h)^{old} - \frac{1}{\alp_m}
c_S^{old}x_1^{S_1} \cdot ... \cdot x_m^{S_m-\alp_m+1} \cdot ...
\cdot x_n^{S_n}F_{x_m}.$$ Hence the change of $c_I$ in this step
is expressed by
\begin{equation} \label{AlgCoChange}
c_I^{new}= c_I^{old}-\frac{I_m-S_m+\alp_m}{\alp_m}c_S^{old}\cdot
a_{I-S+J^m}. \end{equation}
For any basic index $S$, the $R_S$ is equal to the coefficient of
$x^S$ in $p$ after the termination of the algorithm.

Note that $w(I-S+J^m)=w(I)-w(S)+1$, and recall that $a_J=0$ if
$w(J)<1$. Hence $c_I$ is influenced only if $w(S)<w(I)$ and
$w(I-S+J^m)<w(I)$. This proves (i).

Let us now prove by induction that at every step of the algorithm,
$c_I - (1-w(I))b_I \in Bm+G^2$ for any basic coefficient $I$, and
$c_I \in B+G$ for any $I \in E$. After the first step of the
algorithm these statements clearly hold. We suppose that they hold
before a step in which we eliminate $x^S$ using $F_{x_m}$, and
show that they still hold after this step.

First let $I$ be a basic index.

Let us consider several cases.
\begin{enumerate}
\item If $S_m>I_m+1$ then
$a_{I-S+J^m} \in B$ and hence $a_{I-S+J^m} \cdot c_S^{old} \in
Bm$.

\item If $S_m \leq I_m$ then $S$ is a basic index and by induction
hypothesis $c_S^{old} \in B+G^2$ and hence $a_{I-S+J^m} \cdot
c_S^{old} \in Bm+G^2$.

\item If $S_m =I_m +1$ and $I_m < \alp_m-2$ then $S$ is a basic index and by induction
hypothesis $c_S^{old} \in B+G^2$ and hence $a_{I-S+J^m} \cdot
c_S^{old} \in Bm+G^2$.

\item If $S_m =I_m +1$ and $I_m = \alp_m-2$ then $S \in E$ and
$I-S+J^m \in E$ and hence by induction hypothesis $a_{I-S+J^m}
\cdot  c_S^{old} \in (B+G)(B+G) \subset Bm+G^2$.
\end{enumerate}

Now let $I \in E$. Then there exists $k$ such that $I_k =
\alp_k-1$ and $I_j \leq \alp_j - 2$ for $j \neq k$. We know that
$S_m \geq \alp_m -1$, $S_p \leq I_p$ for $p \neq m$ and $I_m-S_m +
\alp_m \geq 0$.

Let us consider several cases.
\begin{enumerate}
\item If $S_m > I_m +1$ then $a_{I-S+J^m} \in B
+G$.
\item If $S_m < I_m +1$ then $c_S^{old}  \in B +G$.
\item If $S_m =I_m +1$ and $S_k >0$ then $a_{I-S+J^m} \in B
+G$.
\item If $S_m =I_m +1$ and $S_k =0$ then $c_S^{old} \in B +G$.
\end{enumerate}

This proves (ii).

Substituting $I=(\alp_1-2,...,\alp_n-2)$ in \ref{AlgCoChange} and
arguing in the same way we obtain (iii). \proofend

\begin{corollary}
The system (\ref{RIzero}) is equivalent to the system consisting
of equations
\begin{equation}\label{SecondEqC1}
b_{I}=\widetilde{\psi}_{I}(g_{J},u_{K})\text{ for }I\neq
(\alp_1-2,...,\alp_n-2),\end{equation}  where
$\widetilde{\psi}_{I}\in G^2$ are polynomials in $g_J$ and
$u_{K}$,

and the last equation
\begin{equation}\label{LastEqC1_N}R:=\sum_{I \in E}
A_k\cdot g_{I}\cdot
 g_{dual(I)}+\Theta(g_{J},u_{K})=0,
\end{equation}
where $A_k$ are positive rational numbers and
$\Theta(g_{J},u_{K})$ is a polynomial from $G^2m$.
\end{corollary}

%


Now we see that our substratum germ is isomorphic to the germ at 0
of the affine variety given in the affine space with coordinates
$\{g_{J},u_{K}\}$ by the last equation (\ref{LastEqC1_N}). The
quadratic part $Q$ of this equation is a non-degenerate quadratic
form in $\{g_{I}|I \in E \}$. We will now show that unless $n=4$
and $\alp_1=\alp_2=\alp_3=\alp_4=3$, the quadratic form $Q$ has
rank at least 3.
   For this it is enough to find 3 points $I\in E$. Indeed for any such point
$I=(j_1,...,j_{k-1},\alp_k-1,j_{k+1},...,j_n)$ the dual point
$I'=(\alp_1-2-j_1,...,\alp_{k-1}-2-j_{k-1},\alp_k-1,\alp_{k+1}-2-j_{k+1},...,\alp_{n}-2-j_n)$
also satisfies these conditions.  Consider several cases
separately:

(0) For $n=4$, $\alp_1=\alp_2=\alp_3=\alp_4=3$, $d=12-9=3$, there
are no $g_{I}$. Moreover, $\mathcal{D}$ is empty. So
$V_{d,H}^{0,0}$ consists of one point. Hence $V_{d,H}$ is also a
smooth algebraic variety of dimension $n^2-1+n=19$. The expected
dimension is $\binom{n+d}{n}-1-\prod_{i=1}^n(\alpha_i-1)=
\binom{7}{4}-1-2^4=18$. So, $V_{d,H}$ has non-expected
codimension. Till
the end of the proof we assume that this is not the case.

(1) $\alp_1\geq 4$. In this case we have the points
$(1,1,0,...,0,\alp_{n}-1)$ and $(2,0,0,...,0,\alp_{n}-1)$ in $I$.
Since $d\geq \alp_1$, $\alp_2\geq 3$. (1.1) $\alp_2=3$. In this
case $n\geq 4$ and $\alp_3\geq 3$ since $d\geq \alp_1$. Hence we
have one more point $(1,0,1,0,...,0,\alp_{n}-1)$  in $I$. (1.2)
$\alp_2\geq 4$. In this case we have point
$(0,2,0,...,0,\alp_{n}-1)$ in $I$.

(2) $\alp_1=3$. In this case $n\geq 4$ and
$\alp_1=\alp_2=\alp_3=\alp_4=3$ for the same reason. So unless
$n=4$ we have points $(1,1,0,...,0,\alp_{n}-1)$,
$(0,1,1,0,...,0,\alp_{n}-1)$ and $(1,0,1,0,...,0,\alp_{n}-1)$ in $I$.\\

Note that for all mentioned points in cases  1 and 2,
$|I|=\alp_n+1$ which is less than or equal to $d$. Note also that
the weight of all these points is at least $\frac{2}{\alp_1}+
\frac{\alp_{n}-1}{\alp_{n}}>
\frac{2}{2\alp_{n}}+\frac{\alp_{n}-1}{\alp_{n}} = 1$.

Since the quadratic form is non-degenerate of rank $\geq 3$, our
substratum germ is a reduced irreducible non-smooth variety of
expected codimension and of order two.

Now we are going to prove that the singular locus $Y$ of our
substratum germ coincides with the germ $X_0$ at $H$ of the affine
subspace $X=Z(G)$. Since equation (\ref{LastEqC1_N}) lies in
$G^2$, all its first order partial derivatives lie in $G$, and
hence the singular locus includes $X_0$. Let $Z$ be the variety
given by the equations
$$\frac{\partial R}{\partial g_{I}}=0$$
for $I\in E$.
 Since the linear part of this system of equations
is non-degenerate, the germ $Z_0$ of $Z$ at $f$ is smooth and
hence irreducible. Clearly $ Y \subset Z_0$ and hence $X_0\subset
Z_0$. They have the same dimension and $Z_0$ is irreducible hence
$X_0=Z_0$ which implies $X_0=Y$.

So our substratum germ is a reduced irreducible non-smooth variety
of expected codimension which has a smooth singular locus.

\subsection{Proof of the theorem for the case $\alp_1\geq
2\alp_n.$} \label{HyperCase2}

We will make now a series of coordinate changes so that in new
coordinates $F$ (see \ref{polinomF}) will not have terms of the
form $x_i^{\alp_i-1}x^J$ lying below and on the Newton polytope,
except $x_i^{\alp_i}$.
The first
one will be $x_i\mapsto x_i$ for $i<n$ and $x_n\mapsto
x_n-\frac{1}{\alpha_n}a_{\langle 2,0,...,0,\alp_{n}-1
\rangle}x_1^{2}$. The coefficients of the polynomial $F$ in the
new coordinates are expressed through the coefficients in the old
coordinates by the formula
\begin{equation} 
a_{\langle i_1,...,i_n\rangle}^{new} =a_{\langle
i_1,...,i_n\rangle}+\sum_{s=1}^{[i_1/2]}(-1)^s \binom{k+s}{s}
a_{\langle i_1-2s,i_2,...,i_{n-1},i_n+s\rangle} (\frac{a_{ \langle
2,0,...,0,\alp_{n}-1 \rangle }}{\alpha_n})^s.
\end{equation}
After this coordinate change the coefficient $a^{new}_{\langle
2,0,...,0,\alp_{n}-1 \rangle}$ will vanish.

Note that in the new coordinates $F$ might get terms of degree
more than $d$. In fact, those terms might have very high degrees.
However, by finite determinacy theorem (Theorem \ref{F-D}) $f$ is
$d+3$-determined. This means that after each coordinate change we
may (and will) erase all the terms of $F$ of degree more than
$d+3$.

In the same way, we get rid of all the coefficients of the form
$a_{\langle j_1,...,j_{n-1},\alp_{n}-1 \rangle}$ in ascending
order of the corresponding monomials. We recall that the monomial
ordering we use is the negative weighted degree lexicographical
ordering with $w=(1/\alp_1,...,1/\alp_n)$ (see Example \ref{mo}).
The coordinate change indexed $J=(j_1,...,j_{n-1},\alp_{n}-1)$
will be $x_i\mapsto x_i$ for $i<n$ and $x_n\mapsto
x_n-\frac{1}{\alpha_n}a_{\langle
j_1,...,j_{n-1},\alp_{n}-1\rangle}x_1^{j_1}\cdot ...\cdot
x_{n-1}^{j_{n-1}}$. The coefficients of the polynomial $F$ in the
new coordinates are expressed through the coefficients in the old
coordinates by the formula
\begin{equation} \label{CoordChange_N2}
a_{\langle i_1,...,i_n\rangle}^{new} =a^{prev}_{\langle
i_1,...,i_n\rangle}+\sum_{s=1}^{min\{[i_l/j_l]|j_l\neq 0, \, 0\leq
l \leq n-1\}}(-1)^s \binom{k+s}{s}a^{prev}_{\langle
i_1-sj_1,...,i_{n-1}-sj_{n-1},i_n+s \rangle}
(\frac{a^{prev}_{J}}{\alpha_n})^s.
\end{equation}
As can be seen from formula (\ref{CoordChange_N2}), $g_{I}$ can be
affected only during a coordinate change whose index does not
exceed $I$ by any coordinate, and hence has lower weight.  Thus
after all these coordinate changes, all coefficients $g_{\langle
j_1,...,j_{n-1},\alp_{n}-1\rangle}$ will be zero. Denote by $a'_I$
the coefficient of the monomial $x^I$ after the coordinate
changes. It is easy to see that $a'_I=a_I+\phi'_I(a_J)$ where
$\phi'_I(a_J)\in mG$. Note that
$$a'_{\langle \alp_1-2,...,\alp_n-2\rangle}= \sum A_J\cdot
g_{\langle j_1,...,j_{n-1},\alp_n-1\rangle}\cdot
 g_{\langle\alp_1-2-j_1,...,\alp_{n-1}-2-j_{n-1},\alp_n-1\rangle}+\Phi'(a_I)$$ where $A_J\in \R\setminus \{0\}$ for all $J$, and
$\Phi'(a_I)\in mG^2$.

We continue with the coordinate changes and, in the same way as
before, we get rid of the coefficients of the form $g'_{\langle
j_1,...,j_{n-2},\alp_{n-1}-1,0\rangle}$ starting from $g'_{\langle
1,0,...,0,\alp_{n-1}-1,0 \rangle}$. This time it might be non-zero
since $\phi'_{\langle 1,0,...,0,\alp_{n-1}-1,0\rangle}(a_J)$ may
be non-zero. However, it lies in the ideal $G$.

We do the same for $g'_{\langle
j_1,...,j_{k-1},\alp_{k}-1,0,...,0\rangle}$ for all $k\geq 2$ in
the descending order of $k$.

Denote by $\wa_I$ the coefficient of the monomial $x^I$ after all
these coordinate changes. Again, $\wa_I=a_I+\phi_I(a_J)$ where
$\phi_I(a_J)\in mG$ and
$$\wa_{\langle\alp_1-2,...,\alp_n-2\rangle}= \sum_{k=1}^n \sum A_J\cdot
g_{\langle j_1,...,j_{k-1},\alp_k-1,0,..,0\rangle}\cdot
 g_{\langle\alp_1-2-j_1,...,\alp_{k-1}-2-j_{k-1},\alp_k-1,\alp_{k+1}-2,...,\alp_n-2\rangle}+\Phi(a_I)$$
where $A_J\neq 0$ for all $J$ and $\Phi(a_J)\in mG^2$.

Now we want to find out for which $\{a_{I}\}$ $H_{F}$ lies in
$V_{d,H}^{0,0}$. Let $\wF(x_1,...,x_n)$ be the polynomial $F$ in
new coordinates. $\wF$ doesn't include monomials $x^I$ for $I$
below and on the Newton polytope and satisfying $I_j= \alp_j-1$
for some $j$. Hence, by Corollary \ref{Var}, in order to belong to
our substratum $\wF$ should include no terms below and on the
Newton polytope except of $x_i^{\alp_i}$.

 In other words, we have the following equations on
$\wa_{I}$:
\begin{align} \label{FirstEqC3_N}&\we_{I}=0 \\
 \label{FirstAdEqC3_N}
&\widetilde{q}_{I}=0
\end{align}

As explained in the previous subsection (and also subsection
\ref{Struct2}), $F \overset{c}{\sim} f$ iff $M_F \simeq T_F$ i.e.
$F(x_1,...,x_n)\in\langle F_{x_1},...,F_{x_n} \rangle$.

As in case one, in order to check that we use the
\textsc{redNFBuchberger} algorithm with the negative weighted
degree lexicographical ordering with $w=(1/\alp_1,...,1/\alp_n)$
(see Example \ref{mo}).

Again, we can stop the algorithm manually when the leading
monomial of the tail is less than $x_1^{\alp_1-2}\cdot...\cdot
x_n^{\alp_n-2}$ (see Algorithm \ref{OurRedNFBuch}).

As a result, we obtain the normal form
$$NF(\wF|\langle \wF_{x_1},...,\wF_{x_n}
\rangle)=\sum R_{I}(\wa_{J})x^I,$$ where $x^I$ are elements of the
basis of algebra $\mathbb{C}[x_1,...,x_n]/_{\langle
f_{x_1},...,f_{x_n} \rangle}$ which lie above $\Delta(f)$ and
 $R_{I}(\wa_{J})$ are polynomials in
$\wa_{J}$. Hence $\wF$ belongs to the ideal $\langle
\wF_{x_1},...,\wF_{x_n} \rangle$ if and only if all the
coefficients $R_{I}(\wa_{J})$ are 0.

Thus we obtain a system of equations on $\wa_{J}$:
\begin{equation}
R_{I}(\wa_{J})=0
\end{equation}
where $R_I$ has the form $R_I= \wb_I
\prod_{i=1}^n\wu_{i}^{\gamma_{I,i}} +
\psi_{I}(\wg_{J},\wu_{K},\wb_L),$ where all $L$ have weight less
than that of $I$. Thus we can, as in case 1, express $\wb_I$ and
obtain an equivalent  system of equations:

\begin{equation}\label{ThirdEqC3_N}
\wb_{I}=\frac{\wpsi_{I}(\wg_{J},\wu_{K})}{\prod_{i=1}^n\wu_{i}^{\gamma_{I,i}}}
\end{equation} where $\wpsi_{I}\in \wG^2$ for
$\wG=\langle \wg_{J}\rangle$ and $\wu_{i}=\wa_{\langle
0,...,0,\alp_i,0,...,0\rangle}$.

The number of equations in system (\ref{ThirdEqC3_N}) is equal to
the number of basis coefficients above the Newton polytope
$\Delta(f)$.

Now we express new coefficients through the old ones. Recall that
$\wu_{i}=1+\phi_{\langle 0,...,0,\alp_i,0,...,0\rangle}(a_{J})$,
$\wa_{I}=a_{I}+\phi_{I}(a_{J})$ for other $I \in \mathcal{D}$, and
$\wa_{I}=\phi_{I}(a_{J})$ for $I\notin \mathcal{D}$, where
$\phi_{I}(a_{J})\in mG$.

Consider  system (\ref{FirstEqC3_N}). After substituting old
coefficients it will be $e_{I}=-\phi_{I}(a_{J})$, where
$\phi_{I}\in m^2$. The same with the system (\ref{FirstAdEqC3_N}).

Consider system (\ref{ThirdEqC3_N}) except the last equation i.e.
the equation on $\wb_{\langle \alp_1-2,...,\alp_n-2\rangle}$.
After substituting the old coefficients and multiplying by
denominators it will become
\begin{equation} \label{Eqs_N}
b_{I}=\widetilde{\psi}_{I}(a_{J}),
\end{equation}
where $\widetilde{\psi}_I \in m^2$.

The last equation is of particular interest. After passing to the
old coordinates and multiplying by the denominator its right hand
side will become
$$\wpsi_{\langle\alp_1-2,\dots, \alp_n-2\rangle}(\wg_{J},\wu_{K})=
\sum_{J \in S } A_J\cdot g_{J}\cdot g_{dual(J)}+\Theta'(a_K),$$
where $S\subset E$ is the set of indexes of terms that did not
vanish during the coordinate changes, $\Theta'(a_{K})$ lies in
$mG^2$ and the $A_J$ are positive rational numbers. The left hand
side will be
$$\wb_{\langle\alp_1-2,\dots,
\alp_n-2\rangle}\cdot\prod_{i=1}^n\wu_{i}^{\gamma_{I,i}}=
(0+\phi_{\langle\alp_1-2,\dots,
\alp_n-2\rangle})(\prod_{i=1}^n(1+\phi_{i})^{\gamma_{I,i}})=\sum_{J
\in T} A_J\cdot g_{J}\cdot g_{dual(J)}+\Psi(a_K),$$ where
$T=\{J\in E|\,J\notin S$, and $dual(J)\notin S\}$, $\Psi(a_K) \in
mG^2$ and $A_J$ are positive rational numbers. Moving the right
hand side to the left we obtain the equation

\begin{equation}\label{LastEqC2} P:=\sum_{J\in E}
B_J\cdot g_{J}\cdot g_{dual(J)}+\Theta(a_{K})=0,
\end{equation}
where $\Theta(a_{K}) \in mG^2$, the $B_J$ are non-zero rational
numbers and $sign(B_J)=sign(B_{dual(J)})$.

So our substratum germ $V_{d,H}^{0,0}$ is isomorphic to the germ
at the origin of the affine variety given in the affine space with
coordinates $\{a_I\}$ by the system of equations

\begin{equation} \label{TotalEqs_N}
e_{I}+\phi_{I}(a_{J})=0
\end{equation}
$$q_{I}+\phi_{I}(a_{J})=0$$
$$ b_{I}-\widetilde{\psi}_{I}(a_{J})=0$$
and the last equation (\ref{LastEqC2}).

We have to prove that $V_{d,H}^{0,0}$ is a reduced irreducible
non-smooth variety of expected codimension which has a smooth
singular locus.

The system  (\ref{TotalEqs_N}) has a diagonal linear part. Hence
all the equations in it are independent and the variety $W$,
defined by it, is smooth at 0. The last equation (\ref{LastEqC2})
does not depend on the preceding ones and does not have a linear
part. Hence $V_{d,H}^{0,0}$ is non-smooth of expected codimension.

The quadratic part of equation (\ref{LastEqC2}) is

$$Q:=\sum_{J\in E} B_J\cdot g_{J}\cdot
 g_{dual(J)}$$

Since the $B_J$ are non-zero, it is a quadratic non-degenerate
form of rank $r$ equal to the number of integer points in $E$.

Now we will show that $r \geq 3$. We do that by exhibiting 3
integer points in $E$: $(2,0,...,0,\alp_n-1)$,
$(1,1,0,...,0,\alp_n-1)$ and $(2,1,0,...,0,\alp_n-1)$. Since
$\alp_1 \leq d$ we have $\alp_2 \geq 3$, and hence $I_j \leq
\alp_j-2$ for $j<n$ for all the 3 points. Hence, it is enough to
show that their degrees don't exceed $d$. Indeed, their maximal
degree is $\alp_n+2 \leq 2\alp_n\leq \alp_1\leq d$. So $r\geq 3$.
Hence the variety defined by the principle part of our system of
equations on $V_{d,H}^{0,0}$ is reduced and irreducible and hence
our germ is reduced and irreducible. Since the quadratic form $Q$
is non-zero $V_{d,H}^{0,0}$ has order two.

Now we are going to prove that the singular locus $Y$ of our
substratum germ coincides with the germ $X_0$ at the origin of the
affine subvariety $X$ given in $W$ by the ideal $G$.

Since  equation (\ref{LastEqC2}) lies in $G^2$, all its first
order partial derivatives lie in $G$, and hence $X_0$ lies in the
singular locus. Consider the jacobian of the system obtained by
merging system (\ref{TotalEqs_N}) with the last equation
(\ref{LastEqC2}). Fix $I\in E$. Let $M(I)$ be the minor of the
jacobian given by columns that include partial derivatives by
$g_I$ and all $e_J$, $b_{J}$ and $q_{J}$. The linear part of
$M(I)$ is $A_Ig_{dual(I)}$. Let $Z_0$ be the germ at the origin of
subvariety of $W$ given by the system of equations $M(I)=0$ for
all $I\in E$. Since the linear part of this system is diagonal,
$Z_0$ is irreducible and has the same dimension as $X_0$. As
$X_0\subset Y\subset Z_0$, it implies $X_0=Y=Z_0$.

So the substratum germ $V^{0,0}_{d,H}$ is reduced irreducible
non-smooth variety of expected codimension which has a smooth
singular locus.
\subsection{Proof that the sectional singularity type is $A_1$ for all cases}

The sectional singularity type of a variety germ with smooth
singular locus is the singularity type of transversal intersection
of the singular locus with a linear space.

Consider the linear subspace $L$ spanned by $e_I$, $q_I$, $b_I$
and $g_I$. We have seen that it is transversal to the singular
locus of $V^{0,0}_{d,H}$ in the affine space with coordinates
$\{a_{I}\}$. Hence the sectional singularity type of
$V^{0,0}_{d,H}$ is the singularity type of the scheme theoretic
intersection $L \cap V^{0,0}_{d,H}$. $L \cap V^{0,0}_{d,H}$ is
given in $L$ by a system of equations (\ref{TotalEqs_N}) with
linear part non-degenerate in $e_I$, $q_I$ and $b_I$ and one more
equation with quadratic principle part which is not degenerate in
the variables $g_I$. Hence the singularity type of $L \cap
V^{0,0}_{d,H}$ is $A_1$. \proofend

\section{Proof of the theorem on stability properties
of obstructed equianalytic families} 
\label{ProofSquares} \setcounter{lemma}{0}
This section is dedicated to the proof of Theorem \ref{Squares}.

\subsection{The structure of the proof}

First of all, we pass to local coordinates $x_i=\frac{t_i}{t_0}$,
and denote $f_0(x_1,...,x_n):=F_0(1,x_1,...,x_n)$.

It is known (Lemma \ref{MathYauCor}) that polynomials
$g_1=f_1+x_{n+1}^2(1+h_1)$ and $g_2=f_2+x_{n+1}^2(1+h_2)$, where
$f_i \in \C\{x_1,...,x_n\}$ and $h_i \in m \subset
\C\{x_1,...,x_{n+1}\}$, are contact equivalent if and only if the
$f_i$ are contact equivalent.

For any polynomial $F = f_0+x_{n+1}^2+\sum a_{I,j}x^Ix_{n+1}^j \in
m^2 \subset \C\{x_1,...,x_{n+1}\}$ there exists an analytic
diffeomorphism of $(\C^{n+1},0)$ that brings $F$ to the form
$f_0+x_{n+1}^2(1+h)$. One can write explicit formulas for this
diffeomorphism that depend polynomially on the coefficients of $F$
(see Section \ref{SqCorChang}). Now we build a map of germs
$\phi:(|\mathcal{O}_{\PP^{n+1}}(d)|,W^1) \to
(|\mathcal{O}_{\PP^n}(\tau+1)|,H)$ in the following way: take the
equation of the hypersurface which includes $x_{n+1}^2$ with
coefficient 1, bring it to the form $f_0+x_{n+1}^2(1+h)$ (using
the above diffeomorphism) and take the hypersurface defined by the
($\tau+1$)-jet of $f_0$. By Lemma \ref{MathYauCor} and the finite
determinacy theorem (Theorem \ref{F-D}), the preimage of
$V_{\tau+1,H}$ will be $V^U_{d,W^1}$.

In Section \ref{SqCorChang} we obtain an explicit formula for
$\phi$ (formula (\ref{NewToOld})). The linear part of $\phi$ is
the identity and the quadratic part depends only on the
coefficients of monomials which include $x_{n+1}$ with degree 1.

Since $V_{\tau+1,H}$ is T-smooth, it is locally defined by a
system (*) of $\tau$ equations having non-degenerate linear part.
Using $\phi$ we obtain a system (**) of equations on
$V^U_{d,W^1}$. This is done in Section \ref{SquareEqs}. Since the
linear part of $\phi$ is the identity, the linear part of (**) is
obtained from the linear part (*) by substituting zeros for the
coefficients of monomials of degree more than $d$. Since the
tangent space to $V_{d,H}$ has codimension $\tau -h^1$, the rank
of the linear part of (**) will be $\tau -h^1$. This proves
statement (a).

The quadratic part of (**) is the sum of two systems. First is
obtained from quadratic part of (*) by substituting zeros for
coefficients of monomials of degree more than $d$. The second
summand is obtained from the linear part of (*) by composing it
with the quadratic part of $\phi$ and then substituting zeros for
the coefficients of the monomials of degree more than $d$. We show
that if $h^1(\kj_{Z^{ea}(H)/\mathbb{P}^{n}}(2d-2))=0$ then the
second summand is non-zero.

In case $h^1=1$, (**) has $\tau-1$ equations with independent
linear parts and one more equation with quadratic principal part
which is independent of the linear parts of the previous
equations. In Section \ref{SqSecPfE} we show that it has rank at
least 2 for $m \geq \max\{1,5-d\}$ and at least 3 for $m \geq
\max\{1,6-d\}$. This finishes the proof of statement (e).

In Section \ref{reduced} we analyze the system (**) in the general
case and prove statement (b).

In Sections \ref{SqSecPfC} and \ref{SqSecPfD} we prove statements
(c) and (d) using Lemmas \ref{Cast} and \ref{Davis} on the
Castelnuovo function.

\subsection{Coordinate changes} \label{SqCorChang}

First of all, we pass to local coordinates $x_i=\frac{t_i}{t_0}$,
and denote $f_0(x_1,...,x_n):=F_0(1,x_1,...,x_n)$.

Let $F(x_1,...,x_{n+1})=f_0+x_{n+1}^2+\sum a_{I,j}x^Ix_{n+1}^j$ be
a polynomial of degree $\leq d$. We want to find equations on the
coefficients $a_{I,j}$ that $F$ should satisfy in order to define
a hypersurface that belongs to $V^U_{d,W^1}$. First of all we may
suppose that the coefficient of monomial $x_{n+1}^2$ is one.

Now we want to get rid of the coefficients of monomials
$x^Ix_{n+1}$. In order to do that we make the following series of
coordinate changes: $x_i\mapsto x_i$ for $1\leq i\leq n$,
$x_{n+1}\mapsto x_{n+1}-1/2a_{I,1}x^I$. After this coordinate
change the new coefficients will be expressed from the previous
ones by the following formula:
\begin{equation}\label{SqCoordChang} a_{J,k}^{(I)}=\sum_{s=0}^{\max \limits
_l\{[J_l/I_l]\}}\binom{s+k}{k}(-\frac{1}{2})^s
a^{prev}_{J-sI,s+k}(a_{I,1}^{prev})^s
\end{equation} We start from $I$ having smallest
degree, and continue in increasing order of degrees. All the
coefficients of the form $a_{J,1}$ influenced during the
coordinate change indexed $I$ have degree more than degree of
$a_{I,1}$ (except $a_{I,1}$ which vanishes). Hence we can continue
making such coordinate changes until $F$ has no coefficients
$a_{I,1}$ of degree less than $\tau+2$. Denote the final
coefficients by $a'_{I,j}$. From the formula (\ref{SqCoordChang})
we see that they can be expressed through the original
coefficients by
\begin{equation}\label{NewToOld}a'_{I,0}=a_{I,0}-1/2 \sum
a_{J,1}a_{I-J,1} - 1/4 a_{I/2,1}^2 + \phi_I(a_{K,k})\end{equation}
where $\phi_I \in \langle a_J \rangle^3$.
\subsection{Equations defining $V^U_{d,W^1}$ and proof of
statement (a)} \label{SquareEqs}
 We want to find equations on the coefficients $a_{I,j}$ that
$F$ should satisfy in order to be analytically equivalent to
$f_0$. By the finite determinacy theorem ( \ref{F-D}) we can
suppose that all coefficients of $F$ of degree more than $\tau +1$
are zero. Then we can present $F$ in the form $$F=f_0+\sum
a'_{I,0}x^I+x_{n+1}^2(1+\sum a'_{I,j}x^Ix_{n+1}^{j-2}).$$

Then, by Lemma \ref{MathYauCor} $F \overset{c}{\sim} f_0$ if and
only if $f_0 \overset{c}{\sim} f_0+\sum a'_{I,0}x^I$. Therefore,
in order for $F$ to lie in $V^U_{\tau+1,W^1}$, $f_0+\sum
a'_{I,0}x^I$ should lie in $V_{\tau+1,H}$. Hence in order to find
the needed equations on $a_{I,j}$ we have to take equations that
define $V_{\tau+1,H}$ inside the linear system
$|\mathcal{O}_{\PP^n}(\tau+1)|$, substitute there $a'_{I,0}$ and
then using formulas (\ref{NewToOld}) to express $a'_{I,0}$ through
the old coefficients $a_{I,j}$.

Let us now investigate the equations on $V_{\tau+1,H}$. Since
$h^1(\kj_{Z^{ea}(H)/\mathbb{P}^{n}}(2d-2))=0$, both $V_{\tau+1,H}$
and $V_{2d-2,H}$ are smooth and have expected codimension $\tau$.
That means that there is a system of $\tau$ local equations on
$V_{\tau+1,H}$ with a non-degenerate linear part which remains
non-degenerate after substituting zeroes in place of all
coefficients of monomials of degree bigger than $2d-2$.

Replacing the system of equations by an equivalent one, we can
suppose that the linear part of this system has echelon form such
that every row has a special element $a_{L_i}$ which appears only
in the linear part of this row and satisfies $|L_i| \leq 2d-2$.

When we substitute in these equations zeros instead of the
coefficients of degree $>d$ the rank of the linear part of the
system drops by $h^1$. That means that $h^1$ rows of the linear
part include only coefficients of degrees $>d$. Renumber the rows
so that those will be the last $h^1$ equations. Denote the index
of the special element of equation number $i$ by $L_i$.

As can be seen from formula (\ref{NewToOld}), the linear part of
$a'_{I,0}$ is $a_{I,0}$, and the quadratic part is $-1/2 \sum
a_{J,1}a_{I-J,1} - 1/4 a_{I/2,1}^2$. Hence the linear parts of the
last $h^1$ equations remain 0 when we express the new coefficients
through the old ones.

Therefore the stratum germ $V^U_{d,W^1}$ is locally defined by a
system of $\tau$ equations, of which only $\tau -h^1$ have a
linear part, which is non-degenerate. Hence
$h^1(\kj_{Z^{ea}(W^1)/\mathbb{P}^{n+1}}(d))=h^1$. By induction the
same is true for all $V^U_{d,W^m}$ for any $m\geq 1$. By Lemma
\ref{MathYauCor} $degZ^{ea}(W^m)=\tau$. Statement (a) is now
proven.

\subsection{Proof of statement (b)} \label{reduced}

We have to show that $V^U_{d,W^m}$ has an irreducible component
which is reduced of expected dimension for $m \geq h^1+1$.
Consider $W^2$. It is obtained from $W^1$ by the same procedure of
adding a square. Hence we know the form of the equations on
$V^U_{d,W^2}$. The system consisting of the first $\tau - h^1$
equations has a non-degenerate linear part, and it is the same
linear part as in the equations on $V_{d,H}$. The last $h^1$
equations start from quadratic forms.

Let us analyze the quadratic form of equation number $\tau -
h^1+j$. Let $w^j_0$ be the quadratic form appearing in the
equation number $\tau - h^1+j$ on $V_{d,H}$ (it may be zero). The
quadratic part of the corresponding equation on $V^U_{d,W^1}$ is
$w^j_0$ + $w^j_1$ where $w^j_1$ depends only on coefficients
$a_{I,1}$, and its form was described above. Hence the quadratic
part of the corresponding equation on $V^U_{d,W^2}$ is $w^j_0$ +
$w^j_1$ + $w^j_2$  where $w^j_1$ depends only on coefficients
$a_{I,1,0}$, and $w^j_2$ is the same as $w^j_1$ but with variables
$a_{I,0,1}$.

Continuing in the same way we see that for general $m$, the
quadratic part of equation number $\tau-h^1+j$ will be
$w^j_0(a_{\langle I,0,\dots,0\rangle})+w^j_1(a_{\langle
I,1,0,\dots,0\rangle})+w^j_2(a_{\langle I,0,1,0,\dots,0\rangle})
+\dots+w^j_m(a_{\langle I,0,\dots,0,1\rangle})$. For convenience,
we denote $0^m:=(0,...,0)\in \Z^m$ and
$e_i:=(0,...,0,1,0,...,0)\in \Z^m$.

Let $X^m$ be the variety defined by the principle part of the
system of equations on $V^U_{d,W^m}$. In order to show that
$V^U_{d,W^m}$ has a reduced component of expected dimension, it is
enough to show that $X^m$ has a reduced component of expected
dimension. To do that we will prove that there is a point in $X^m$
in which the jacobian of the principle part of the system of
equations on $V^U_{d,W^m}$ has maximal rank.

For every $1 \leq j \leq h^1$ we choose $J_j$ and $K_j$ (not
necessary different) such that $J_j+K_j=L_{\tau-h^1+j}$. Here,
$L_i$ is the index of the special element of the linear part of
equation number $i$ on $V_{\tau,H}$. Let $A$ be the linear
subspace spanned by  all $\{a_{K_j,e_{j}}\}_{j=1}^{h^1}$ and by
all $\{a_{I,e_{m}}\}_{1\leq|I|\leq d-1}$.

Consider the minor defined by derivatives
with respect to $a_{L_i,0^m}$ for $1\leq i \leq \tau-h^1$ and to
$a_{J_j,e_{j}}$ for $1\leq j \leq h^1$. 
We claim that the restriction of this minor on the linear subspace
$A$ is $C \cdot  a_{K_1,e_1} \cdot \dots \cdot
a_{K_{h^1},e_{h^1}}$ where $C$ is non-zero real number. This minor
is a determinant of a block matrix in which the first block is the
identity matrix $Id_{(\tau-h^1) \times (\tau-h^1)}$. It is left to
show that the second block is a diagonal matrix with $C_j \cdot
a_{K_j,e_j}$ on the diagonal.

Indeed, consider, for example, equation number $\tau-h^1+1$. The
derivative of $w^1_0$ w.r.t. $a_{J_1,e_1}$ is 0 since $w^1_0$
doesn't depend on it at all. The same is true about $w^1_{\geq
2}$. The derivative of $w^1_1$ w.r.t. $a_{J_1,e_1}$ contains
$a_{K_1,e_1}$ with non-zero coefficient $C_1$. It may also contain
other $a_{I,e_1}$, but they vanish on our subspace $A$. Consider
now the derivative of equation number $\tau -h^1 +j $ (for $j>1$)
w.r.t. $a_{J_1,e_1}$. Again, the derivatives of $w_0$ and
$w^{j}_{\geq 2}$ is zero. The derivative of $w^j_1$ does not
contain $a_{K_1,e_1}$ since the coefficient $a_{J_1+K_1}=a_{L_1}$
does not appear in the linear part of equation number $\tau-h^1+j$
on $V_{\tau,H}$. Hence this derivative also vanishes on the
subspace $A$ spanned by all $\{a_{K_j,e_{j}}\}_{j=1}^{h^1}$ and by
all $\{a_{I,e_{m}}\}_{1\leq|I|\leq d-1}$.

By the same reason, the restriction to $A$ of the derivative of
equation number $\tau -h^1 +j$ with respect to $a_{J_i,e_i}$ is
equal to 0 if $j \neq i$ and $C_i \cdot a_{K_i,e_{i}}$ if $j=i$.

Hence on $A$ the minor is a determinant of the diagonal matrix
with entries $(C_ja_{K_j,e_{j}})$ and hence the minor is
$Ca_{K_1,e_1} \cdot \dots \cdot a_{K_{h^1},e_{h^1}}$. Clearly,
every neighborhood of 0 in $X^m \cap A$ contains a point in which
$Ca_{K_1,e_1} \cdot \dots \cdot a_{K_{h^1},e_{h^1}} \neq 0$. At
this point, a maximal minor of the jacobian is non-zero, hence the
jacobian has maximal rank. Therefore the tangent space to $X^m$ at
this point has codimension $\tau$ and hence $X^m$ has a reduced
component of expected codimension $\tau$. Hence $V^U_{d,W^m}$ also
has a reduced component of expected codimension for $m \geq
h^1+1$.

\subsection{Proof of statement (e)} \label{SqSecPfE}

Now we have $h^1(\kj_{Z^{ea}(H)/\mathbb{P}^{n}}(d))=1$. By Lemma
\ref{Cast}, (d) on the Castelnuovo function
$h^1(\kj_{Z^{ea}(H)/\mathbb{P}^{n}}(d+1))=0$. Hence from Section
\ref{SquareEqs} $V^U_{d,W^1}$ is defined by a system of $\tau-1$
equations with non-degenerate linear part and one more equation
without linear part (since $h^1=1$). Hence $V^U_{d,W^1}$ is
non-smooth.

In order to show that $V^U_{d,W^m}$ has expected codimension
(respectively is reduced, respectively irreducible) it is enough
to show that the scheme defined by the principle parts of the
above system of equations  has expected codimension (respectively
is reduced, respectively irreducible).

For $1 \leq i \leq \tau-1$, we can just express $a_{L_i}$ from
equation number $i$. Therefore, this scheme is isomorphic to the
scheme $X^m$ defined by the quadratic part $w$ of the last
equation in the affine space of coefficients $\{a_{I,J} | 2\leq
|I+J| \leq d\}\setminus \{a_{L_i,0}\}_{i=1}^{\tau -1}$. Since $w$
is non-zero, $X^m$ has expected codimension. $X^m$ is reduced if
$rank(w)\geq 2$ and irreducible if $rank(w)\geq 3$.

As in Section \ref{reduced}, $w=w_0(a_{\langle
I,0,\dots,0\rangle})+w_1(a_{\langle
I,1,0,\dots,0\rangle})+w_2(a_{\langle I,0,1,0,\dots,0\rangle})
+\dots+w_m( a_{\langle I,0,\dots,0,1\rangle})$ where
$rank(w_1)=\dots=rank(w_m)\geq 1$.
Hence $X^m$ is reduced for $m\geq 2$ and irreducible for $m\geq
3$. It is left to deal with $m=1$ and $m=2$. It is enough to show
that for $d= 4$ we have $rank (w_1) \geq 2$ and for $d\geq 5$ we
have $rank (w_1)\geq 3$.

We start with the case $d\geq 5$. The quadratic form $w_1$ can be
expressed by $$w_1= \sum_{(I,J)\, s.t. \, |I+J|=d+1} C_{I,J}
a_{I,e_1}a_{J,e_1}.$$ Note that if $I+J=L_{\tau}$ then $C_{I,J}
\neq 0$. Since $|L_{\tau}|=d+1$, it can be presented as
$L_{\tau}=I_1+I_2$ where $|I_1|=2$ and $|I_2|=d-1$ and also as
$L_{\tau}=I_3+I_4$, where $|I_3|=3$ and $|I_4|=d-2$. Note that
$I_1$, $I_2$ and $I_3$ are different since they have different
degrees since $d \geq 5$. If $d=5$, it is possible that $I_3=I_4$.
Consider the reduction $w_1'$ of $w_1$ on the (3 or 4 dimensional)
linear subspace spanned by
$a_{I_1,e_1},a_{I_2,e_1},a_{I_3,e_1},a_{I_4,e_1}$. The only sums
of pairs of those multiindexes whose degrees are $d+1$ are
$I_1+I_2$ and $I_3+I_4$. Hence
$w_1'=C_{L_{\tau}}(a_{I_1,e_1}a_{I_2,e_1}+a_{I_3,e_1}a_{I_4,e_1})$
which has rank 3 or 4. Hence the rank of $w_1$ is at least 3.
Hence $V^U_{d,W^m}$ is reduced and irreducible for $m\geq 1$.

Consider now $d = 4$. As in the previous case, we can find $I_1$
and $I_2$ such that $|I_1|=2$, $|I_2|=d-1=3$ and
$I_1+I_2=L_{\tau}$. By reducing $w_1$ on the subspace spanned by
$a_{I_1,e_1}$ and $a_{I_2,e_1}$ we see that the rank of $w_1$ is
at least 2. Hence $V^U_{d,W^m}$ is reduced for $m\geq 1$ and
irreducible for $m\geq 2$.

\subsection{Proof of statement (c)} \label{SqSecPfC}

By Lemma \ref{Cast} (d), it is enough to show that
$\mathcal{C}_{Z^{ea}(H)}(2d-3))=0$. Let $Z(j(f_0))$ be
zero-dimensional scheme defined by the jacobian $j(f_0)$ where
$f_0$ is the local equation of $H$. Then $Z^{ea}(H)$ is a
subscheme of $Z(j(f_0))$ and hence by Lemma \ref{Cast} (e) it is
enough to show that $\mathcal{C}_{Z(j(f_0))}(2d-3))=0$. Let
$C_1=Z(\frac{\partial f_0}{\partial x_1})$ and
$C_2=Z(\frac{\partial f_0}{\partial x_2})$ and let $Z'$ be the
complete intersection $C_1 \cap C_2$. Then $Z(j(f_0))$ is a
subscheme of $Z'$ and hence by Lemma \ref{Cast} (e) it is enough
to show that $\mathcal{C}_{Z'}(2d-3))=0$. Let $k$ be the degree of
$\frac{\partial f_0}{\partial x_1}$ and $l$ be the degree of
$\frac{\partial f_0}{\partial x_2}$. By Lemma \ref{Davis},
$\mathcal{C}_{Z'}(l+k-1))= 1-1=0$. Since $k,l \leq d-1$ we obtain
$\mathcal{C}_{Z'}(2d-3)=0$. Therefore
$h^1(\kj_{Z^{ea}(H)/\mathbb{P}^{2}}(2d-4))=0$ and hence by (b) the
germs $V^U_{d,W^m}$ have a reduced component of expected dimension
for $m\geq h^1+1$.

\subsection{Proof of statement (d)} \label{SqSecPfD} Suppose
$h^1(\kj_{Z^{ea}(H)/\mathbb{P}^{n}}(2d-2))> 0$. Then
$h^1(\kj_{Z^{ea}(H)/\mathbb{P}^{n}}(l))> 0$ for all $d\leq l \leq
2d-2$. Hence by Lemma \ref{Cast},(d) the Castelnuovo function
$C_{Z^{ea}(H)}(l)>0$. Hence
$h^1(\kj_{Z^{ea}(H)/\mathbb{P}^{n}}(d))-h^1(\kj_{Z^{ea}(H)/\mathbb{P}^{n}}(2d-2))=\sum_{l=d+1}^{2d-2}C_{Z^{ea}(H)}(l)\geq
d-2$. Hence $h^1(\kj_{Z^{ea}(H)/\mathbb{P}^{n}}(2d-2))\leq
h^1(\kj_{Z^{ea}(H)/\mathbb{P}^{n}}(d))-(d-2) \leq 0$. So
$h^1(\kj_{Z^{ea}(H)/\mathbb{P}^{n}}(2d-2))=0$ and by (b) the germs
$V^U_{d,W^m}$ have a reduced component of expected dimension for
$m\geq h^1+1$. \proofend

\begin{remark}
{\rm If one wishes $W^m$ to be unisingular, they can be defined by
$F_0 +\sum_{j=1}^m t_{n+j}^2t_0^{d-2}+\sum_{j=1}^m
\lambda_jt_{n+j}^d$, where the $\lambda_j$ are generic. The same
proof shows that an analogous theorem will hold about the
equisingular family germ $V_{d,W^m}$. }
\end{remark}

\begin{remark}
{\rm Statement (c) can be strengthened as follows: let $H$ be a
projective plane curve of degree $d$ which is not a union of $d$
lines through the same point. Then
$h^1(\kj_{Z^{ea}(H)/\mathbb{P}^{2}}(2d-5))=0.$}
\end{remark}
For a proof see \cite{Gou2}, Remark 3.4.4.


\end{document}